\documentclass[11pt]{amsart} 
\usepackage[all]{xy}
\theoremstyle{remark}
\newtheorem{rem}{Remark}

\theoremstyle{definition}

\theoremstyle{theorem}
\newtheorem{theo}{Theorem}[section]
\newtheorem{cor}[theo]{Corollary}
\newtheorem{prop}[theo]{Proposition}
\newtheorem{lem}[theo]{Lemma}

\newcommand{\bib}[6]{ \bibitem{#1}  #2, \textit{ #3\/}, #4 \textbf{#5} #6}

\begin{document}

\title[Connecting Manifolds and Hopf Invariants]{Homotopical
Dynamics II: Hopf Invariants, Smoothings and the Morse Complex.}
\author{Octavian Cornea}
\address{Universit\'{e} de Lille 1\\
U.F.R. de Math\'{e}matiques \& U.R.A. 751 au CNRS\\
59655 Villeneuve D'Ascq, France}

\email{cornea@gat.univ-lille1.fr}
\urladdr{http://www-gat.univ-lille1.fr/\~{}cornea/octav.html}

\date{December 10, 1998}

\begin{abstract}
The ambient framed bordism class of the connecting manifold of two consecutive
critical points of a Morse-Smale function is estimated by means of a certain Hopf
invariant. Applications include new examples of non-smoothable Poincar\'e duality
spaces as well as an extension of the Morse complex.
\end{abstract}

\maketitle

\section{Introduction.}
Let $M$ be a smooth, compact, riemannian manifold and let 
$f:M\longrightarrow {\bf R}$ be a smooth Morse-Smale function, regular and
constant on $\partial M$. The flow $\gamma:M\times {\bf R}\longrightarrow
M$ used below is induced by $-\nabla f$. Assume that
$P$ and $Q$ are consecutive critical points of
$f$ (this means that $f(P)> f(Q)$ and that there are no broken flow lines connecting
$P$ to $Q$) of indexes, respectively,
$p$ and $q$.

\

An important, classical problem in Morse theory is to use the topology of $M$ to
understand the properties of the moduli space $Z(P,Q)$ of flow lines 
that connect $P$ to $Q$. 

\

This problem is the main motivation of the paper. The key new idea introduced here
 is that knowledge of the homotopy of the based loop space $\Omega M$ leads to  
significant information on these moduli spaces.

\

More precisely, recall that the genericity of the
Morse-Smale condition implies that $Z(P,Q)$ is a manifold of dimension
$p-q-1$ called connecting manifold of $P$ and $Q$. It has a canonical
normal framing and a classical result of John Franks \cite{Franks} claims that its
framed bordism class $\{Z(P,Q)\}\in \Omega^{fr}_{\ast}$ 
 is given (via the Thom-Pontryagin construction) by the relative attaching map
$\delta_{f}(P,Q)$ associated to the succesive  attachments of
the cells corresponding to the critical points $Q$  and $P$. 

The closure of the
space of all the points situated on some flow line joining 
$P$ to $Q$ is identified to the unreduced suspension $\Sigma Z(P,Q)$. Therefore,
we have an inclusion $\Sigma Z(P,Q)\hookrightarrow M$ and if $M$ is
simply connected (which we will asssume) a well defined adjoint
$$l(P,Q):Z(P,Q)\longrightarrow \Omega M$$
 The above mentioned normal framing of
$Z$ together with the map $l(P,Q)$ provide, via the Thom-Pontryagin construction,
a homotopy class $T(P,Q):S^{p-1}\longrightarrow \Sigma^{q}(\Omega M ^{+})=
S^{q}\vee S^{q}\wedge \Omega M$. By Franks' result the projection of $T(P,Q)$
on $S^{q}$ is $\delta_{f}(P,Q)$. Let $h(P,Q)\in \pi_{p-1}
(\Sigma^{q}\Omega M)$ be the projection of $T(P,Q)$ on the second factor.

\

The main result of the paper continues the work of Franks by giving a purely
homotopical, computable description of the ambient, framed, bordism class
of $Z(P,Q)$, $[Z(P,Q)]^{fr}\in\Omega_{p-q-1}^{fr}(\Omega M)$.  
In fact, we show that $\Sigma h(P,Q)$
equals the suspension of a certain Hopf invariant, $H(P,Q)$, associated to the
succesive cell attachments corresponding to $Q$ and
$P$.  As  $T(P,Q)$ represents $[Z(P,Q)]^{fr}$ we conclude that this bordism class
equals the stable image of $H(P,Q)+\delta_{f}(P,Q)$. 

\

The homotopy classes $h(P,Q)$ turn out to be highly relevant for understanding the
respective connecting manifolds and also for the topology of $M$ itself. We
study some of their properties. 

Specializing to the Morse-Smale
case the Spanier-Whitehead duality results of \cite{Cornea1}, \cite{Cornea2}, we
show that, stably, $\epsilon\delta_{-f}(Q,P)$ equals $\delta_{f}(P,Q)+\Delta (P,Q)$
where
$\epsilon\in\{-1,+1\}$ with the twisting term
$\Delta (P,Q)=J^{q}(\Omega^{q}\Sigma^{q}(\Omega \nu)\circ h^{\ast}(P,Q))$. Here
$\nu:M\longrightarrow {\bf BSO}$ classifies the stable normal bundle of $M$,
$h^{\ast}(P,Q)$ is the $q$-th order adjoint of $h(P,Q)$ and
$J^{q}:\pi_{k}(\Omega^{q}\Sigma^{q}\mathbf{SO})\longrightarrow
\pi_{k}^{S}$ is a factor of the classical J-homomorphism
$J:\pi_{k}(\mathbf{SO})\longrightarrow \pi_{k}^{S}$. This can be used to
deduce a measure of the embedding complexity of $Z(P,Q)$.
We also use this result to construct examples of non-smoothable Poincar\'e duality
spaces (many of which are PL-manifolds). When the obstructions
 to smoothing concern only the relative attaching maps of the top cell the relevant
morphism is, of course, $J=J^{0}$ and one recovers in this case classical
results (see for instance \cite{Spivak}) and the examples of Smith \cite{Smith} 
which are thus seen to be part of a more general pattern. 

In a different direction, if $P$, $R$ are critical points such that $f(P)>f(R)$ and all
 flow lines connecting $P$ to $R$ are broken at most once, let $I(P,R)$
be the set of intermediate critical points. We show that the relation
$\sum_{Q\in I(P,R)}(-1)^{(p+r)q}[Z(P,Q)]^{fr}\bullet [Z(Q,R)]^{fr}=0$
 is satisfied inside the ring $\Omega^{fr}_{\ast}(\Omega M)$ (where
the product $\bullet$ is induced by loop composition).
 This can be used to define various chain complexes which, in the simplest possible
case (when the index of successive critical points differs by just one), provide the
Morse complex of $f$. 

\

The material is organized as follows. After a second section, containing recalls and
notations, comes the technical heart of the paper, in section three, very much in
the spirit of the work of Franks \cite{Franks}. The description of $\Sigma h(P,Q)$
in terms of the Hopf invariant is given at this point, as well as the proof of the relations
among the bordism classes of the connecting manifolds. The fourth section contains
applications and most homotopy theoretical arguments are concentrated here.  

\

{\bf Acknowledgements.} I thank Raoul Bott, Fred Cohen, Pascal Lambrechts, Chuck
McGibbon and Alberto Verjovsky for useful discussions as well as John Harper
whom I also thank  for pointing out \cite{Smith}. I am most grateful to Mark
Mahowald for his valuable suggestions and for his encouragement. 

\section{Recalls and notations.}

\subsection{Hopf invariants.}

For further use we fix some standard facts. We denote by $A\ast B$ the join of $A$
and $B$. For pointed spaces, $A\ast B\simeq \Sigma A\wedge B$, $\Sigma
A=S^{1}\wedge A$. Recall also that $A^{+}$ is the pointed
space obtained by the disjoint union of $A$ and a disjoint base point; if $A$ is
pointed we have $(A\times B)/(\ast\times B)=A\wedge B^{+}$. 
For $A$ and $B$ connected we will also need the existence of the fibration
$\Omega A\ast\Omega B\longrightarrow A\vee B\longrightarrow A\times B$
which is trivial after looping.

Assume that $X$ is a $CW$-complex and that
$X'\hookrightarrow X''\stackrel{i}{\hookrightarrow} X$ are connected subcomplexes
such that  there is a cofibration
sequence $S^{q-1}\stackrel{f}{\longrightarrow}X'\longrightarrow X''$. Fix also a
map $S^{p-1}\stackrel{g}{\longrightarrow}X''$.

The Hopf invariants that we will use are defined using this
data (they are versions of invariants defined by Ganea \cite{Ganea}\cite{Dula}).

The Hopf invariant of $g$ relative to $f$ is obtained as follows.
Consider the map $$t:S^{p-2}\longrightarrow \Omega
S^{p-1}\stackrel{\Omega g}{\longrightarrow}\Omega X''
\stackrel{\Omega\nabla}{\longrightarrow}\Omega
(S^{q}\vee X'')\stackrel{p}{\longrightarrow}\Omega (\Omega
S^{q}\ast\Omega X'')$$ Here, $\nabla$ is the coaction $X''\longrightarrow
S^{q}\vee X''$ and $p$ is the canonical projection in the splitting
$\Omega(S^{q} \vee X'')\simeq
 \Omega S^{q}\times \Omega X'' \times \Omega (\Omega S^{q}\ast\Omega X'')$; 
the first map in the composition is the restriction to the bottom cell.

We have a projection $r:\Omega S^{q}\ast\Omega X''
\longrightarrow S^{q}\wedge \Omega X''$ induced by the canonical evaluation
$\Sigma\Omega S^{q}\longrightarrow S^{q}$.

\

\textit{The needed Hopf invariant is the homotopy class:
$$H(g,f)=(id_{S^{q}}\wedge \Omega i)\circ r\circ t^{\ast}: S^{p-1}\longrightarrow
S^{q}\wedge \Omega X$$  where 
$t^{\ast}:S^{p-1}\longrightarrow \Omega S^{q}\ast \Omega X''$ is the
adjoint of $t$}. 

We denote by $\delta=\delta (f,g) : S^{p-1}\longrightarrow S^{q}$ the relative
attaching map given by projecting $\nabla\circ g$ onto $S^{q}$. 

\

\begin{rem} For example, let $\ast : S^{1}\longrightarrow \ast$ and let
$\eta:S^{3}\longrightarrow S^{2}$ be the Hopf map and $X=\mathbf{CP}^{2}$.
Then $H(\eta,\ast):S^{3}\longrightarrow S^{2}\wedge\Omega\mathbf{CP}^{2} $
is the inclusion of the bottom cell.  Similarly,
let $\ast : S^{q-1}\longrightarrow S^{t}$ be the trivial map and
$w : S^{q+t-1}\longrightarrow S^{q}\vee S^{t}$ be the obvious Whitehead
product and $X=S^{q}\times S^{t}$. Then $H(w,\ast) : S^{q+t-1}\longrightarrow
 S^{q}\wedge\Omega (S^{q}\times S^{t})$ is the inclusion
$S^{q}\wedge S^{t-1}\longrightarrow  S^{q}\wedge\Omega S^{t}
\longrightarrow  S^{q}\wedge ( \Omega S^{q}\times \Omega S^{t})$.
\end{rem} 

\subsection{Elements of Morse theory.} 
We follow here the fundamental paper of John Franks \cite{Franks}.

Let $M^{n}$ be a smooth compact manifold and let $f : M\longrightarrow
\mathbf{R}$ be a smooth function. If $\partial M\not=\emptyset$ we assume that the
function is constant and regular on $\partial M$. We assume also that a
riemannian metric is fixed on $M$ and we denote by $\gamma : M\times
\mathbf{R}\longrightarrow M$ the flow induced by $-\nabla f$. A critical point $P$ of $M$ is non-degenerate
if $Hess_{P}(f)$ is a non-degenerate matrix. The index of the induced bilinear
form is called the index of $P$. The set $W^{u}(P)=\{x\in M : 
\lim_{t\rightarrow -\infty}\gamma_{t}(x)=P\}$ is called the 
{\em unstable manifold} of
$P$ and $W^{s}(P)=\{x\in M : \lim_{t\rightarrow +\infty}\gamma_{t}(x)=P\}$
is the {\em stable manifold} of $P$. If $P$ is non-degenerate and of index $p$,
then $W^{s}(P)\approx Int(D^{n-p})$ and $W^{u}(P)\approx Int(D^{p})$. 

We assume
from now on that $f$ is Morse which means that all its critical points
are non-degenerate, and even Morse-Smale which means that if $P$ and $Q$
are two critical points of $f$, then $W^{u}(P)$ and $W^{s}(Q)$ are in general
position. The Morse-Smale condition is generic. 

We say that two critical points $P$ and $Q$ are {\em consecutive} if 
$f(P)> f(Q)$ and there are no broken flow lines connecting $P$ to $Q$.
 In this case, let $f(Q)< a <f(P)$ and let $S^{s}(Q)=W^{s}(Q)\bigcap f^{-1}(a)$ and
$S^{u}(P)=W^{u}(P)\bigcap f^{-1}(a)$. If $P$ and $Q$ are critical and 
consecutive we may (and will) assume also (possibly after
 slightly isotoping the function $f$) that $P$ and $Q$ are the only critical points in
 $f^{-1}([f(Q),f(P)])$. In this case, it is easy to see that $S^{u}(P)\approx
S^{p-1}$ and $S^{s}(Q)\approx S^{n-q-1}$.
The Morse-Smale condition insures that $S^{u}(P)$ and
$S^{s}(Q)$ intersect transversely.  Their intersection, $Z(P,Q)$, is called the
connecting manifold of $P$ and $Q$ (a.k.a. the moduli space of
flow lines connecting $P$ to $Q$).
It is a $p-q-1$ -dimensional manifold. Fix for each critical point $R$ of
$f$ an orientation on the linear subspaces $V^{u}(R)$ and $V^{s}(R)$ of
$T_{R}(M)$ that are respectively tangent to $W^{u}(R)$ and $W^{s}(R)$. If the
manifold $M$ is oriented, pick these orientations such that at each point $R$ they give
on $V^{u}(R)\oplus V^{s}(R)$ the fixed orientation of $T_{R}(M)$. It is easy to see
that a choice of a basis of $V^{u}(R)$ induces a framing of the normal bundle of
$W^{s}(R)$  and, similarly, a choice of basis for $V^{s}(R)$ induces a normal framing of
$W^{u}(R)$. We choose these bases in a way compatible with the fixed orientations
and then these framings are unique up to isomorphism.  We now return to the
two consecutive critical points $P$ and $Q$. The normal bundle of
$Z(P,Q)$ in $S^{u}(P)$ is induced by the normal bundle of $W^{s}(Q)$
in $M$ and hence inherits a standard framing (coming from the choices in
$V^{u}(Q)$). 

It is well known that the passage through a non-degenerate critical point corresponds
to the attachment of a cell of dimension the index of the critical point.
In other words, if $M'=f^{-1}(-\infty, f(Q)-\epsilon]$,
$M''=f^{-1}(-\infty,f(P)-\epsilon]$, $M'''=f^{-1}(-\infty, f(P)+\epsilon]$ with
$\epsilon$ small enough we have cofibration sequences:
$S^{q-1}\stackrel{\alpha(Q)}{\longrightarrow} M'\longrightarrow M''$
and $S^{p-1}\stackrel{\alpha(P)}{\longrightarrow} M''\longrightarrow M'''$.
Denote by $\delta_{f}(P,Q) : S^{p-1}\longrightarrow S^{q}$ the obvious relative
attaching map. 

One of the key results in \cite{Franks} is
that {\em $\delta_{f}(P,Q)$ corresponds to the standard framing of $Z(P,Q)$ in
$S^{u}(P)$ via the Thom-Pontryagin construction}. When
$p-q=1$ this comes down to just counting (with sign)
 the number of elements in $Z(P,Q)$. Hence, as immediate application of
this result of Franks, one obtains that if $C_{i}={\bf Z}/2 < x : \nabla f (x)=0,
ind(x)=i>$ are ${\bf Z}/2$ vector spaces and $d:C_{i}\longrightarrow C_{i-1}$ is the
unique linear application defined on basis elements by
$d(x)=\sum_{ind(y)=i-1}\#(Z(x,y))y$, then $(C_{i},d)$ is a complex, called the
Morse complex of $f$, and 
$H_{\ast}(C_{i},d)\simeq H_{\ast}(M;{\bf Z}/2)$ (in the oriented case, by
using some appropriate signs when counting the elements of $Z(x,y)$ in the definition
of $d$, one obtains the integral homology of $M$). There are at least a couple of
other proofs of this fact. An analytical one appears in
\cite{Witten}. Another, that will be extended in sections 3 and 4,  is based on 
understanding the boundaries of certain moduli spaces of
connecting flow lines (see for example
\cite{Schwarz}) .

For a smooth Morse-Smale function $f$ on $M$ 
 there is an associated $CW$-decomposition of $M$ 
(we consider now only functions that are constant, regular and maximal on $\partial
M$). If $f$ is continuously deformed via Morse-Smale functions to a second
Morse-Smale function $f'$, then the $CW$-decompositions associated to $f$ and $f'$
are equivalent. When the metric is allowed to vary, the
$CW$-decomposition corresponding to $f$ is determined up to a 
contractible choice \cite{Klein}.

Conversely, if $M$ and $\partial M$ are simply-connected, $n>5$ and
$H_{\ast}(M;{\bf Z})$ torsion free, then all $CW$-complexes $Y$ of the homotopy
type of $M$ and satisfying a certain minimality condition correspond to some
Morse-Smale function on $M$ \cite{Franks} (the minimality condition is the
following: 
$Y$ has a unique $0$-dimensional cell, at most one $n$- cell and
if $e^{k}_{1}$ and $e^{l}_{2}$ are two cells
of $Y$ with $k>l$, then either $2k-l<n-1$ or the smallest subcomplex of
$Y$ that contains $e^{k}_{1}$ contains also $e^{l}_{2}$). In particular,
there are self-indexed Morse functions having a number of critical
points of index $k$ equal to the rank of $H_{k}(M ;{\bf Z})$ \cite {Smale}. They
are called perfect Morse functions. 

\subsection{Duality and Flows.}
Let $X$ be a $CW$-complex and consider an orthogonal fiber bundle of 
rank $n$ over $X$, $\mu:E\longrightarrow X$. Let $E_{S}(X)\longrightarrow X$
be the associated spherical bundle and $T^{\mu}(X)$ the respective Thom space.
Assume that for a certain cell-decomposition of $X$ we have cofibration
sequences $S^{q-1}\longrightarrow X'\longrightarrow X''$ and
$S^{p-1}\longrightarrow X''\longrightarrow X'''$ with $X'''$ a subcomplex of $X$.
Denote by $\delta : S^{p-1}\longrightarrow S^{q}$ the respective relative
attaching map. It is easy to see that the $CW$-decomposition of $X$
induces one for $T^{\mu}(X)$. For example, by pulling back the
cofibration sequence $S^{q-1}\longrightarrow X'\longrightarrow X''$ to the
spherical fibration of $\mu$ one obtains a push-out square:
$$\xy
\xymatrix@+0cm{
S^{q-1}\times S^{n-1}\ar[d]\ar[r]&E_{S}(X')\ar[d]\\
S^{n-1}\ar[r]&E_{S}(X'') }
\endxy $$

by pushing down this square into the original cofibration we get the
push-out square:
$$\xy
\xymatrix@+0cm{
S^{q-1}\wedge S^{n}\vee S^{n}\ar[d]\ar[r]& T^{\mu}(X')\ar[d]\\
S^{n}\ar[r]&T^{\mu}(X'') }
\endxy $$

As the left vertical map is just the projection onto $S^{n}$
it is easy to transform this push-out square into a cofibration sequence
$S^{n+q-1}\longrightarrow T^{\mu}(X')\longrightarrow T^{\mu}(X'')$.

In particular we obtain a relative attaching map $\delta^{\mu}:S^{n+p-1}
\longrightarrow S^{n+q}$.

We return now to the context and notations of the previous sub-section.
Thus $f:M\longrightarrow \mathbf{R}$ is a smooth Morse-Smale function;
$P$ and $Q$ are consecutive critical points of $f$ of indexes respectively $p$ and
$q$; $\delta_{f}(P,Q)$ is the corresponding relative attaching map. The function $-f$ is
also Morse-Smale and $Q$ and $P$ are consecutive critical points for $-f$. Thus, we
also have a relative attaching map $\delta_{-f}(Q,P) : S^{n-q-1}\longrightarrow
S^{n-p}$.

Let $\nu$ be the stable normal bundle of $M$. The results in \cite{Cornea1}
imply that {\em $\delta_{f}(P,Q)^{\nu}$ agrees stably up to sign with
$\delta_{-f}(Q,P)$ and if $\nu$ is trivial,  
then $\delta_{f}(P,Q)^{\nu}$ and $\delta_{f}(P,Q)$
agree stably}.

\begin{rem} 
When $\nu$ is
trivial the result appears already in the paper of Franks \cite{Franks}.
In fact, one has much more general results valid for general flows and
isolated invariant sets (in the context of Conley index theory) in which the equality up
to sign is replaced by Spanier-Whitehead duality \cite{Cornea1}. 
\end{rem}

\subsection{The {\bf J}-homomorphism.} We need to recall a few elements of
classical homotopy theory. We denote by $\pi^{S}_{k}$ the stable $k$-stem. 
Similarly, the stable homotopy groups of any space $X$ are denoted
by $\pi_{k}^{S}(X)$.

The  J-homomorphism
$J:\pi_{k}({\bf SO})\longrightarrow \pi_{k}^{S}$ is defined as follows. For
$\alpha\in\pi_{k}({\bf SO})$ there is some $m\in{\bf N}$ such that 
$\alpha\in \pi_{k}({\bf SO}(m))$. Consider the composition $$\alpha':S^{m-1}\times
S^{k}\stackrel{id\times \alpha}{\longrightarrow}S^{m-1}\times {\bf
SO}(m)\stackrel{\cdot}{\longrightarrow} S^{m-1}$$

Apply the Hopf construction to this map $\alpha'$ (or, equivalently, suspend
$\alpha'$ and use the splitting of the suspension of the domain to restrict to
$S^{m+k}$) thus getting a map $\alpha'':S^{m+k}\longrightarrow S^{m}$. The
image of $\alpha''$ in the stable $k$-stem is $J(\alpha)$.

One can show \cite{Whitehead} that this construction does not depend of the different
choices involved and that the resulting map is a group homomorphism. 

As the image of  J is stable, it is immediate to see that this homomorphism
factors as $\pi_{k}({\bf SO})\longrightarrow \pi_{k}^{S}({\bf SO})\stackrel{J'}
{\longrightarrow}\pi_{k}^{S}$. The homomorphism $J'$ is called
the "very stable $J$" and also the "bi-stable J-homomorphism"
\cite{Knapp}. It is easy to see that, for each $q\geq 0$, in between $J$ and $J'$ there
is an intermediate factor $J^{q}:\pi_{k}(\Omega^{q}\Sigma^{q}{\bf
SO})\longrightarrow \pi_{k}^{S}$. These are also homomorphisms and they
commute with the morphisms induced in homotopy by the inclusions
$\Omega^{q}\Sigma^{q}{\bf SO}\longrightarrow \Omega^{q+1}\Sigma^{q+1}
{\bf SO}$; $J^{0}=J$ and $J'$ is the limit of the $J^{q}$'s.

For later use we give an explicit description of $J^{q}$. Let $\alpha\in
\pi_{k}(\Omega^{q}\Sigma^{q}{\bf SO})$. There is some $m\in {\bf N}$
such that $\alpha\in \pi_{k}(\Omega^{k}\Sigma^{k}({\bf SO}(m)))$.
Then $J^{q}(\alpha)$ is the stable image of the composition
$$S^{k+m+q}\stackrel{\Sigma^{m}(\alpha^{\ast})}{\longrightarrow}
\Sigma^{m+q}{\bf SO}(m)\stackrel{\mu'}{\longrightarrow}S^{m+q}$$
where $\alpha^{\ast}$ is the $q$-th order adjoint of $\alpha$ and
$\mu'$ is the $q$-th suspension of the Hopf construction applied to the
multiplication $\mu:S^{m-1}\times {\bf SO}(m)\longrightarrow S^{m-1}$.

\begin{rem} The image of  J has been computed by Adams \cite{Adams},
Quillen \cite{Quillen} and Sullivan \cite{Sullivan}.
By the Kahn-Priddy theorem \cite{KahnP} it is known that $J'$ is surjective at the
prime $2$. However, $J'$ is not surjective at any odd prime \cite{Knapp}.
\end{rem}

\section{Morse-theoretic interpretation of the Hopf invariants}
 
As before, let $M^{n}$ be a smooth, riemannian, compact manifold and $f :
M\longrightarrow \mathbf{R}$ a smooth Morse-Smale function regular, maximal
 and constant on $\partial M$.

We assume from now on that $M$ is {\em simply connected} and that
$f$ has a single local minimum.

Suppose that $P$ and $Q$ are consecutive critical points  of
$f$ (in particular $f(P)> f(Q)$) of indexes, respectively, $p$ and $q$. 

For $q\geq 1$, as in 2.2, consider the induced
cofibration sequences:
$S^{q-1}\stackrel{\alpha(Q)}{\longrightarrow} M'\longrightarrow M''$
and $S^{p-1}\stackrel{\alpha(P)}{\longrightarrow} M''\longrightarrow M'''$.
In this case, denote by $H(P,Q)$ the Hopf invariant $H(\alpha(P),\alpha(Q)):
S^{p-1}\longrightarrow \Sigma^{q}\Omega M$. 
When $q=0$ the relevant cofibration sequence is $S^{p-1}\stackrel{\alpha(P)}
{\longrightarrow}\ast=M''\longrightarrow S^{p}=M'''\hookrightarrow M$
and we let $H(P,Q)$ be the adjoint of the inclusion $S^{p}\hookrightarrow M$.

Notice also that $K(P,Q)=\overline{W^{u}(P)\bigcap W^{s}(Q)}$ is homeomorphic
to the (un)reduced  suspension $\Sigma Z(P,Q)$.  As $M$ is simply connected there is
a canonical adjunct (up to homotopy) of the inclusion $i_{Z} : \Sigma Z(P,Q)
\hookrightarrow M$ that we denote by $l(P,Q) :Z(P,Q)\longrightarrow \Omega M$.
By making the choices described in 2.2 we have a standard framing of the normal
bundle of $Z(P,Q)$ in $S^{p-1}\approx S^{u}(P)$. 
Consider an inclusion
of a tubular neighborhood $U\approx D^{q}\times Z\subset S^{p-1}$  induced by this
framing. We fix an orientation of $Z$ such that this inclusion is oriented.  The
Thom-Pontryagin construction applied to the fixed framing together with the map
$l(P,Q)$ gives the map $T(P,Q):S^{p-1}\longrightarrow
S^{p-1}/\overline{S^{p-1}-U}=(D^{q}\times Z)/(S^{q-1}\times
Z)\stackrel{\Delta}{\longrightarrow} (D^{q}\times Z)/(S^{q-1}\times
Z)\wedge Z^{+}\stackrel{p_{1}\wedge l(P,Q)}{\longrightarrow} S^{q}\wedge
(\Omega M )^{+}=S^{q}\vee S^{q}\wedge \Omega M$ ($\Delta$ is induced
by the diagonal and $p_{1}$ is the projection on $S^{q}$). Recall that
$h(P,Q)$ is the projection of $T(P,Q)$ on $S^{q}\wedge \Omega M$
and that the projection of $T(P,Q)$ on $S^{q}$ is $\delta_{f}(P,Q)$.  

\begin{theo} We have the equality: $\Sigma h(P,Q)=\Sigma
H(P,Q)$. In particular, $[Z(P,Q)]^{fr}\in\Omega_{p-q-1}^{fr}(\Omega M)$
equals the stable image of $H(P,Q)+\delta_{f}(P,Q)$.
\end{theo} 

\begin{proof}
When $q=0$ the statement is immediate as $K(P,Q)=S^{p}\hookrightarrow
M$. Assume from now on $q>0$. If $p=q+1$ we have $h(P,Q)=H(P,Q)=0$
hence assume also that $p>q+1$. Denote by $Z_{i}(P,Q)$ the connected
components of $Z(P,Q)$ and notice that $h(P,Q)$ is also given by the
composition $S^{p-1}\stackrel{t}{\longrightarrow}\vee_{i}S^{q}\wedge Z_{i}(P,Q)
\stackrel{\Sigma^{q}l(P,Q)}{\longrightarrow}S^{q}\wedge \Omega M$.
Here $t$ is a degree one map given by the sum of the maps
$S^{p-1}\longrightarrow (D^{q}\times Z_{i}(P,Q))/(S^{q-1}\times Z_{i}(P,Q))
=S^{q}\vee S^{q}\wedge Z_{i}(P,Q)
\stackrel{p_{2}}{\longrightarrow}S^{q}\wedge
Z_{i}$.

The proof has five steps.

\subsubsection*{Some special neighborhoods of a critical point.}
As above, assume that
for
$\epsilon >0$ sufficiently small, 
$P$ and
$Q$ are the only critical points in $f^{-1}([f(Q)-\epsilon, f(P)+\epsilon])$. As before,
let $M'=f^{-1}(-\infty, f(Q)-\epsilon]$ and $M''=f^{-1}(-\infty, f(P)-\epsilon]$.
Recall that $\gamma$ is the flow induced by $-\nabla f$ and let $D(x,r)$ be the
closed disk in $M$ of radius $r$ and center $x$.
Let $U'_{\tau,\epsilon}=\{x\in M : f(Q)-\epsilon\leq f(x)\leq
f(Q)+\epsilon , \exists t\in\mathbf{R}\bigcup \{+\infty, -\infty\}$
 such that $ \gamma_{t}(x)\in D(Q,\tau)\bigcap f^{-1}(f(Q))\}$.

It is useful to recall at this time that, by the Morse lemma, the topology of $f$ inside
$U'_{\tau,\epsilon}$ is independent of $\tau$, for $\tau$ sufficiently small.

 Let $U'=U'_{\tau,\epsilon}$ for a small, fixed $\tau$. 
Then $U'\approx
D^{q}\times D^{n-q}$, $\partial U'= A'\bigcup B'\bigcup C'$ with
$S^{q-1}\times D^{n-q}\approx A'=U'\bigcap f^{-1}(f(Q)-\epsilon)$, $D^{q}\times
S^{n-q-1}\approx C'= U'\bigcap f^{-1}(f(Q)+\epsilon)$, $B'\approx S^{q-1}\times
S^{n-q-1}\times [0,1]$ and $\nabla f$ is tangent to $B'$ in all points $x\in Int (B')$ (we
consider the $0$-end of $B'$ to be contained in $A'$) . Of course, $M''\simeq
M'\bigcup_{A'} U'$. It is obvious that there is a second neighborhood $U\subset U'$ very
close to $U'$, homeomorphic to $U'$, whose boundary admits the same structure as that of
$U'$ and such that if we denote by
$A$, $B$, $C$ the respective pieces of the boundary of $U$, then $A=A'$,
$C\subset C'$, $C\approx C'$, $B\approx B'$ and $\nabla f$ is transverse to $B$ and
points inside $U$.

\subsubsection*{The attaching map $\alpha (P)$.}We may define a deformation
retract $r : M''\times [0,1]\longrightarrow M''$ whose $1$-end $r_{1}$ sends each
point in $M''$ along the flow $\gamma$ (induced by $-\nabla f$) to the point 
where it first  reaches $M'\bigcup_{A}U$. Let us denote $Z=Z(P,Q)$. Consider the
inclusion $S^{p-1}\subset M''$ that represents the relative attaching map $\alpha (P)$.
We see that $S^{p-1}\bigcap C\approx D^{q}\times Z$ (where we choose the
framing of the respective tubular neighborhood of $Z$ in $S^{p-1}$ as described in
2.2). With this identification, we have $\partial D^{q}\times Z=S^{q-1}\times
Z\subset \partial C\approx S^{q-1}\times S^{n-q-1}$. By using the deformation $r$ 
we may assume that the image of $\alpha (P)$ lies in $M'\bigcup_{A} U$ and 
that $\alpha (P) (D^{q}\times Z, S^{q-1}\times Z)\subset (C,\partial C)$. Moreover,
we may consider a collar neighborhood $V$ of $S^{q-1}\times Z$ in $S^{p-1}- Int
(D^{q}\times Z)$ such that $V\approx S^{q-1}\times Z\times [0,1]$ with
$S^{q-1}\times Z\times \{1\}\subset \partial D^{q}\times Z$ and with the
property that (with the identifications described above)  
$\alpha (P)|_{V}:S^{q-1}\times Z\times [0,1]\longrightarrow B$ equals
 $\alpha (P)|_{(S^{q-1}\times Z)}\times
id _{[0,1]}$  and $\alpha (P) (S^{p-1}-(V\bigcup D^{q}\times Z))\subset M'-A$.

Denote by $V'=V\bigcup D^{q}\times Z$, $M^{*}=M'\bigcup_{A} U$ and let
$T=\overline{S^{p-1}-V'}$.   Of course, $V'\approx D^{q}\times Z$.

We now intend to describe the composition $$c: S^{p-1}\stackrel{\alpha
(P)}{\longrightarrow}M''\stackrel{\nabla}{\longrightarrow}S^{q}\vee M''
\longrightarrow  S^{q}\vee M$$ 
Here, $\nabla$ is the coaction and the last map
is induced by inclusion. By the definition of $\nabla$ and making use of $r$
this map is homotopic to $S^{p-1}\stackrel{\alpha (P)}{\longrightarrow}
M''\stackrel{r_{1}}{\longrightarrow} M^{*}\stackrel{\nabla'}{\longrightarrow}
M^{*}/A\stackrel{w}{\longrightarrow} S^{q}\vee M^{*}\longrightarrow S^{q}\vee
M$. Here $\nabla'$ is the obvious collapsing map and the last map is, as before,
the inclusion.  The homotopy equivalence $w$ is the inverse of the obvious one
obtained from the fact that $A\simeq S^{q-1}$ and using the standard flow induced
nullhomotopy of $A\hookrightarrow M^{*}$ (this is defined by first collapsing $A$ to
its core
$S^{q-1}= S^{u}(Q)$ and then collapsing this one to $Q$ along flow lines).  In other
words we have $w: M^{*}/A\stackrel{id}{\longrightarrow} U/A\vee
M'/A\stackrel{w'}{\longrightarrow}S^{q}\vee M^{*}$ whith $w'$ respecting the
wedge, its restriction to $U/A$ is 
$U/A=(D^{q}\times D^{n-q})/(S^{q-1}\times
D^{n-q})\stackrel{p_{1}}{\longrightarrow} D^{q}/S^{q-1}= S^{q}$.
In the wedge $S^{q}\vee M^{\ast}$ the point $Q\in M^{\ast}$ is identified
with the image of $S^{q-1}$ in the quotient $D^{q}/S^{q-1}$.
The restriction of $w'$ to $M'$ sends $A$ to $Q$ by means of a map
$w'':M'\longrightarrow  M'\bigcup W^{u}(Q)\hookrightarrow M^{*}$ which is the
inclusion outside a neighborhood of $A$, is defined inside this neighborhood by using
the null-homotopy mentioned above and induces a homotopy equivalence
$M'/A\longrightarrow M^{*}$. 

\subsubsection*{Description of $c$.} The map $c$ can be described by writting
$S^{p-1}=V'\bigcup _{\partial V'} T$ and giving its restrictions to
each of these two pieces.  On $V'$ the map $c$ is induced by  the map of pairs
$(V'\approx D^{q}\times Z, S^{q-1}\times Z)\hookrightarrow (U, A)$ and composition
with $w'$; this sends $\partial T=\partial V'=S^{q-1}\times Z$ to $Q$. On $T$
it is defined by $k=w''\circ r_{1}\circ \alpha (P)$. 

We will now see that the composition
$k':T\stackrel{k}{\longrightarrow}M^{*}\hookrightarrow M$ 
is homotopic $rel (\partial T)$ to a map $k''$ defined as follows.
Consider the inclusion of pairs $(D^{p},S^{p-1})\hookrightarrow (M,M'')$
whose restriction to $S^{p-1}$ is $\alpha (P)$ (here $D^{p}\subset W^{u}(P)$,
$D^{p}\bigcap M''=\partial D^{p}$).  There is a deformation
$d:D^{p}\longrightarrow M$ of this inclusion that is induced by the flow
$\gamma$, collapses $Z\subset S^{p-1}$ to $Q$ along $\gamma$, is constant ouside a
neighborhood of $Z$ in $D^{p}$ and factors as $D^{p}\longrightarrow D^{p}\bigcup
W^{s}(Q)
\hookrightarrow M$. As $D^{p}\bigcup (W^{s}(Q)\bigcap M'')\simeq \Sigma Z$ that
means that, up to homotopy, $d$ factors through $i_{Z}$.  In $S^{p-1}$ there is
a deformation $l:T\longrightarrow S^{p-1}$ that
is constant outside a neighborhood of $\partial T$ and that sends the point
$(x,y)\in S^{q-1}\times Z=\partial T$ to $(0,y)\in Z$. The map
$k''$ is given by $i_{Z}\circ d\circ l$. To see that $k'$ and $k''$ are homotopic,
notice that they are both homotopic ($rel$ boundary) to the following map
$k'''$: transport $T$ homeomorphically along the flow till it reaches $f^{-1}(f(Q))$.
Let $T'$ be the image of $T$ inside this singular hypersurface.
Use the conical structure of $f^{-1}(f(Q))$ around $Q$ to deform $\partial T'$ to $Q$
without leaving $f^{-1}(f(Q))$ and whithout moving points that are outside a small
neighborhood of $\partial T'$. 

It follows that we may use the map $k''$ instead of $k'$ in the description of $c$. 
We get the following commutative diagram:
$$ \xy
\xymatrix@-0.5cm{
S^{q-1}\times Z\ar[rr] \ar[dd] \ar[dr] & & D^{q}\times Z \ar@{.>}[dd]\ar[dr] &\\ 
& T \ar[dd]^{l} \ar[rr]& & S^{p-1} \ar[dd]^{\alpha'(p)} \\ 
S^{q-1}\times Z\ar@{.>}[rr] \ar[dd] \ar[dr]
& & D^{q}\times Z\ar@{.>}[dd]\ar@{.>}[dr] &\\
 &  D^{p}\ar[dd]^{d} \ar[rr]& & \vee_{i}(S^{q}\vee S^{q}\wedge Z_{i})
\ar[dd]^{j}\\ Q\ar@{.>}[rr] \ar[dd] \ar[dr] & & S^{q} \ar@{.>}[dd] \ar@{.>}[dr] &\\
 &  \Sigma Z \ar[dd]^{i_{Z}} \ar[rr]& & S^{q}\vee \Sigma Z \ar[dd] \\
Q\ar@{.>}[rr] \ar[dr] & &S^{q}\ar@{.>}[dr] &  \\
&  M\ar[rr]& &S^{q}\vee M}
\endxy $$

In this diagram all horizontal squares are push outs and the vertical maps in the right
corner are induced by the other three; $Z_{i}$ are the connected components of $Z$.
We have identified $V'$ to
$D^{q}\times Z$. The composition originating in $T$ is the map $k''$ described above
and that defined on  $D^{q}\times Z$ is the restriction of $c$ (which, as mentioned
above is the projection on $D^{q}$ followed by the collapsing to
$D^{q}/S^{q-1}=S^{q}$). Therefore, the composition in the right corner is
homotopic to $c$. The map $\alpha'(P)$ is additive with respect to the connected
components of $Z$. Its projection onto $\vee_{i}\Sigma^{q}Z_{i}$
 is  homotopic to the degree one map $t$. Thus, the Hopf invariant being
additive (in the first variable), we may assume from now on $Z$ connected.

\subsubsection*{Identification of a Whitehead product.}
With this assumption, the next step is to consider the map
$j: S^{q}\vee S^{q}\wedge Z
\longrightarrow  S^{q}\vee\Sigma Z$ of the diagram above and show that its
restriction to $S^{q}\wedge Z$ is homotopic to the (generalized) Whitehead product
of the inclusions $S^{q}\hookrightarrow  S^{q}\vee\Sigma Z$ and $\Sigma
Z\hookrightarrow S^{q}\vee \Sigma Z$. For this consider the
next commutative diagram.
$$ \xy
\xymatrix@-0.5cm{
S^{q-1}\times Z\ar[rr] \ar[dd] \ar[dr] & & D^{q}\times Z \ar@{.>}[dd]\ar[dr] &\\
 & S^{q-1}\times CZ \ar[dd] \ar[rr]& & S^{q}\wedge Z \ar[dd] \\ 
S^{q-1}\times Z\ar@{.>}[rr] \ar[dd] \ar[dr] & &
D^{q}\times Z\ar@{.>}[dd]\ar@{.>}[dr] &\\
 &  D^{p}\ar[dd] \ar[rr]& &  S^{q}\vee S^{q}\wedge Z \ar[dd] \\
Q\ar@{.>}[rr] \ar[dr] & & S^{q}  \ar@{.>}[dr] &\\
 &  \Sigma Z \ar[rr]& &  S^{q}\vee\Sigma Z}
\endxy $$

Again, the horizontal squares are push outs and the vertical maps in the right
corner are induced by the respective three others, $CZ$ is the cone on $Z$.
The map $S^{q-1}\times Z\longrightarrow D^{p}$ which is the restriction of
$l$ to $\partial T$ factors as $S^{q-1}\times Z\subset S^{q-1}\times CZ
\stackrel{p_{2}}{\longrightarrow} CZ\hookrightarrow D^{p}$ (the last
map being the inclusion of $CZ$ in $W^{u}(P)$, the vertex of this cone being
identified to $P$) and it is this factorization that is used in the upper, vertical, left
square. It is clear that the top, vertical map in the right corner is, up to homotopy, the
inclusion on the first factor. Now, the composition $S^{q-1}\times CZ\longrightarrow
D^{p}\longrightarrow
\Sigma Z$ is given by projection onto $CZ$ and collapsing onto $\Sigma Z$. Similarly,
$D^{q}\times Z\longrightarrow S^{q}$ is projection onto $D^{q}$ and
then collapsing onto $S^{q}$. This shows, by the definition of the 
Whitehead product \cite{Whitehead}, that the composition in the right corner is the
wanted Whitehead product.

We therefore obtain that $c$ is homotopic to the composition
$$c' : S^{p-1}\stackrel {\alpha'(P)}{\longrightarrow} S^{q}\vee S^{q}\wedge Z
\stackrel{h}{\longrightarrow}  S^{q}\vee\Sigma Z\stackrel{id\vee
i_{Z}}{\longrightarrow}S^{q}\vee M$$
 with $h= id_{S^{q}}\vee [i_{S^{q}}, i_{\Sigma Z}]$. 

\subsubsection*{Identification of the Hopf invariant.}
The last step of the proof is to use the factorization of
$c'$ to evaluate the relevant Hopf invariant. 
We now look to $\Omega (h)$ and
use the standard splitting of the loop space of a wedge to write the Hopf invariant
$H(P,Q)$ as the adjoint of the composition
$h': S^{p-2}\stackrel{i}{\hookrightarrow} \Omega
S^{p-1}\stackrel{u}{\longrightarrow}
\Omega  (S^{q}\wedge Z)
\times \Omega (\Omega S^{q}\ast\Omega (S^{q}\wedge
Z))\stackrel{v}{\longrightarrow}
\Omega (\Omega S^{q}\ast\Omega\Sigma
Z)\stackrel{y}{\longrightarrow}\Omega (S^{q}\wedge\Omega\Sigma Z )
\longrightarrow \Omega ( S^{q}\wedge\Omega M)$.

By the basic properties of the Whitehead product we
obtain that the restriction of $h''=y\circ v$ to $\Omega (S^{q}\wedge Z)$ is just the
looping of the map $S^{q}\wedge Z\longrightarrow 
S^{q}\wedge \Omega\Sigma Z$ induced by
the inclusion $Z\hookrightarrow\Omega\Sigma Z$. This shows that the Hopf
invariant verifies $H(P,Q)=\phi +\phi'$. Here,
$\phi :S^{p-1}\stackrel{t}{\longrightarrow} S^{q}\wedge Z \longrightarrow
 S^{q}\wedge\Omega\Sigma Z \longrightarrow  S^{q}\wedge\Omega M$.
As the composition $S^{q}\wedge Z\longrightarrow 
 S^{q}\wedge\Omega\Sigma Z \longrightarrow  S^{q}\wedge\Omega M$ is the
$q$-th suspension of the adjunction of the inclusion $i_{Z}$ we obtain that
$\phi =h(P,Q)$.

The second homotopy class, $\phi'$, is the
composition $S^{p-1}\stackrel{s}{\longrightarrow}
\Omega S^{q}\ast \Omega (S^{q}\wedge Z)\stackrel{z}{\longrightarrow}
\Omega S^{q}\ast\Omega\Sigma Z\longrightarrow  S^{q}\wedge\Omega M$ 
with $s$ being the adjunction of
$p_{2}\circ u\circ i$ and $z$ the top composition in the next diagram.
$$ \xy
\xymatrix@-0.3cm{
\Omega S^{q}\ast\Omega (S^{q}\wedge Z)\ar[r]^{x}\ar[d]&
\Omega(S^{q}_{1}\vee S^{q}_{2})\ast\Omega\Sigma
Z\ar[r]\ar[d]&
\Omega S^{q}\ast\Omega\Sigma Z\ar[d]\\
 S^{q}_{1}\vee S^{q}\wedge Z\ar[d]\ar[r]^{id\vee [,] }&
S^{q}_{1}\vee S^{q}_{2}\vee\Sigma Z\ar[d]\ar[r]^{m\vee id}&
S^{q}\vee\Sigma Z\ar[d]\\  S^{q}_{1}\times S^{q}\wedge Z\ar[r]&
(S^{q}_{1}\vee S^{q}_{2})\times\Sigma Z
\ar[r]&
 S^{q}\times\Sigma Z}
\endxy $$

Here, the columns are fibrations and the top row is induced by the bottom two;
$m$ is the folding map of the two spheres $S^{q}_{1}$ and $S^{q}_{2}$;
the bottom, left, horizontal map is trivial when projected on $\Sigma Z$. 

We now compose $z$ with the evaluation $ev: \Sigma\Omega
S^{q}\wedge \Omega\Sigma Z\longrightarrow S^{q}\wedge \Omega\Sigma Z$.  
The composition $ev\circ z$ factors via the evaluation $ev': \Sigma\Omega
(S^{q}_{1}\vee S^{q}_{2})\wedge
\Omega\Sigma Z\longrightarrow  (S^{q}_{1}\vee S^{q}_{2})\wedge\Omega\Sigma
Z= (S^{q}_{1}\wedge \Omega\Sigma Z)\vee (S^{q}_{2}\wedge\Omega\Sigma Z)$.
Notice that the projection of $ev'\circ x$ on each of the factors of the wedge is null.
This implies that $\Sigma \phi'=0$ and concludes the proof.
\end{proof}

\begin{rem} It is instructive to consider the case of the 
Morse-Smale function $f:S^{2}\times S^{2}\longrightarrow {\bf R}$
with precisely four critical points. Denote by $P$ the maximum and let 
$Q$ be one of critical points of index two.  The equality
$\Sigma H(P,Q)=\Sigma h(P,Q)$ comes down to the fact that
the Thom-Pontryagin construction applied disjointly to
two circles embedded with linking number
one and trivially framed in $S^{3}$ produces the Whitehead product
$S^{3}\longrightarrow S^{2}\vee S^{2}$. 
\end{rem}

Here is a first context in which $h(P,Q)$ is relevant.
Suppose $f$, $P$, $Q$ are as in the theorem and fix a fiber bundle on $M$ 
that is classified by a map $\mu : M\longrightarrow \mathbf{BSO}(m)$ with
$p-2q<m$.  For $x\in \Omega^{fr}_{\ast}(\Omega M)$ let $\overline{\mu}(x)\in 
\Omega^{fr}_{\ast}(\Omega M)$ be defined as follows. Consider
$X\stackrel{g}{\longrightarrow}\Omega M$ together with a framing
$X\stackrel{i}{\hookrightarrow}S^{k}$ representing $x$. The class
$\overline{\mu}(x)$ is represented by $X\stackrel{g}{\longrightarrow}\Omega M$
together with a framing
$X\stackrel{i}{\hookrightarrow}S^{k}\stackrel{j}{\hookrightarrow}S^{k+m}$
given at a point $a\in X$ by $((\Omega\mu\circ g)(a)(i^{\ast}(j)_{a}),i)$
(here $j$ is the standard framing of $S^{k}$ in $S^{k+m}$).

 In a tubular neighborhood $Y$  of $M$
in the total space of $\mu$ we may consider a function $g:Y\longrightarrow \mathbf{R}$
giving the square of the distance from $M$. If $Y$ is sufficiently small this function
is non-degenerate in the direction of the fibre and, in particular, the
difference $f'=f\circ p-g$ is Morse (here $p:Y\longrightarrow M$ is the restriction
of the projection of the bundle). The critical points $P$ and $Q$ are again consecutive
nondegenerate critical points of $f'$. Their indexes are respectively $m+p$ and
$m+q$. Let $Z^{\mu}(P,Q)$ be the (framed) connecting manifold of $f'$
(as a space it coincides with $Z(P,Q)$ but its framing might be different). Recall that
$\delta_{f}(P,Q)^{\mu}$ is the relative attaching map induced on the Thom space of
$\mu$ by $\delta_{f}(P,Q)$. It is easy to see \cite{Cornea1} that
$\delta_{f'}(P,Q)=\delta_{f}(P,Q)^{\mu}$.

For $x\in \pi_{k+n}(S^{n})$ we denote by
$\{x\}\in \pi_{k}^{S}$ its stable image.

\begin{prop} In $\Omega^{fr}_{p-q-1}(\Omega M)$ we have
$[Z^{\mu}(P,Q)]^{fr}=\overline{\mu}([Z(P,Q)]^{fr})$. Moreover, if
$h^{\ast}(P,Q)\in \pi_{p-q-1}\Omega^{q}\Sigma^{q}\Omega M$ denotes the $q$-th
order  adjoint of $h(P,Q)$, then the stable difference
$\{\delta_{f}(P,Q)^{\mu}\}-\{\delta_{f}(P,Q)\}$ equals $$J^{q}(
\Omega^{q}\Sigma^{q}(\Omega\mu)\circ h^{\ast}(P,Q))$$
\end{prop} 

\begin{proof} We use the notations fixed in the proof of the theorem. As seen above,
the bordism class of $Z(P,Q)$ in
$\Omega M$ is given by a map
$S^{p-1}\stackrel{\alpha}{\longrightarrow}S^{q}\wedge Z^{+}
\stackrel{id\wedge l(P,Q)}{\longrightarrow}S^{q}\wedge (\Omega M)^{+}$. 
Of course, the bordism class of $Z^{\mu}(P,Q)$ is given by a similar
map $S^{m+p-1}\stackrel{\alpha'}{\longrightarrow}S^{m+q}\wedge Z^{+}
\stackrel{id\wedge l(P,Q)}{\longrightarrow}S^{m+q}\wedge (\Omega M)^{+}$.
In general, the map $\alpha'$ is not the suspension of $\alpha$. 
It is immediate that, as in \cite{Cornea2},  
$\alpha'=e\circ\Sigma^{m}\alpha$ where $e:
 S^{m+q}\wedge Z^{+}\longrightarrow 
 S^{m+q}\wedge Z^{+}$ is a homotopy equivalence that appears at the
passage (along the flow $\gamma_{1}$ induced by $-\nabla f'$ ) from a neighborhood
of $P$ to one of $Q$. 
More precisely $e$ is induced by the map of pairs 
$e': (D^{m}\times D^{q}\times Z, S^{m-1}\times D^{q}\times Z)
\longrightarrow (D^{m}\times D^{q}\times Z, S^{m-1}\times D^{q}\times Z)$
that takes $(x,y,z)$ to $(\mu^{\ast}(z)(x),y,z)$ where $\mu^{\ast}: Z
\stackrel{l(P,Q)}{\longrightarrow}\Omega M
\stackrel{\Omega\mu}{\longrightarrow} \mathbf{SO}(m)$. The map $e'$ induces
$e$ by collapsing $S^{m+q-1}\times Z$ to a point. The first part of the statement
is now clear.  Moreover, notice that
$p_{1}\circ e$ restricted to the $q+m$ skeleton (which is a wedge of
$q+m$-dimensional spheres in number equal to the number of connected components
of $Z$) is the identity because the bundle is oriented.  Its restriction to
$S^{q+m}\wedge Z$ is $\Sigma^{q-1}J(\mu^{\ast})$. Here, 
$J(\mu^{\ast}):S^{m}\wedge Z\longrightarrow S^{m}$ is defined by the Hopf
construction on the map $S^{m-1}\times Z\longrightarrow S^{m-1}$ given by
$(x,y)\longrightarrow \mu^{\ast}(y)(x)$. Because $p-2q<m$, 
$\delta_{f'}(P,Q)=\Sigma^{m}\delta_{f}(P,Q)+\Sigma^{q-1}J(\mu^{\ast})\circ t$.
The statement follows from the definition of $J^{q}$.
\end{proof}

\begin{rem} a. It is clear that, as $\Sigma h(P,Q)=\Sigma H(P,Q)$, we may replace
in the formula above $h(P,Q)$ by $H(P,Q)$. When $p>q+1$ notice also that
$[h(P,Q)]= [H(P,Q]$ when viewed in $H_{p-q-1}(\Omega M)$ . 
Moreover, in this case, $[h(P,Q)]$ is the fundamental
class of $Z(P,Q)$ in $\Omega M$.

b. The second part of the above proposition, with $H(P,Q)$ in the place
of $h(P,Q)$,  is also a consequence of purely homotopical results of Dula \cite{Dula}. 

c. The proof of the proposition is in fact
the direct specialization to the Morse-Smale case of a result established in
\cite{Cornea2} for reasonable critical points (a class that contains all isolated, 
analytic singularities).
\end{rem}

For the next result we consider two critical points $P$ and $R$ of $f$ (which is
a function as before) such that $P$ and $R$ are not necessarily consecutive but
$f(P)>f(R)$ and if $Q$ is a critical point such that $Q\in \overline{W^{u}(P)}\bigcap
\overline{W^{s}(R)}$, then $P$ and $Q$ are consecutive and so are $Q$ and $R$
(in other words all broken conecting flow lines between $P$ and $R$ are
broken in just one point). We denote by $I(P,R)$ the set of all such intermediate 
critical points $Q$ associated to $P$ and $R$. 

For $P$ and $Q$ consecutive of indexes respectively $p$ and $q$ recall that we
denote by $[Z(P,Q)]^{fr}\in \Omega_{p-q-1}^{fr}(\Omega M)$ the ambient
bordism class of $Z(P,Q)$. We assume $M$ oriented and make the choice of
orientations described in 2.2. As before, assume the connecting manifolds $Z(P,Q)$
oriented such that the standard framing inside
$S^{p-1}$  summed with this orientation gives the standard orientation of
$S^{p-1}$.

\begin{theo} For the choice of orientations described above we have:
$$\sum_{Q\in I(P,R)}(-1)^{(p+r)q} [Z(P,Q)]^{fr}\bullet [Z(Q,R)]^{fr}=0$$
where $\bullet$ is the product in $\Omega_{\ast}^{fr}(\Omega M)$.
\end{theo}

\begin{proof} 
We may assume whithout loss of generality that $I(P,R)$ is contained is the
same critical level $f^{-1}(c)$. Let $$K(P,R)=
\{x\in M : \lim_{t\rightarrow+\infty}(\gamma_{t}(x))=R,  \lim_{t\rightarrow -\infty}
(\gamma_{t}(x))=P\}\bigcup \{P,R\}$$

The proof is based on understanding how the ends of $K(P,R)$ are embedded in
$M$.

\subsubsection*{Identification of the ends of $K(P,R)$.} Let $p$, $r$, $q$ be
repectively the indexes of
$P$, $R$ and $Q\in I(P,R)$. 
Consider $a\in {\bf R}$ such that $f(R)<a<f(P)$ and let 
$Z(P,R)=f^{-1}(a)\bigcap K(P,R)$. Then $Z(P,R)$ is a
manifold of dimension
$p-r-1$ (of course, $Z(P,R)$ is not closed) and its homeomorphism type does not
depend on the choice of $a$. Clearly, $K(P,R)$ is homeomorphic to the unreduced
suspension  of $Z(P,R)$. We also recall the notation
$K(P,Q)=\overline{W^{u}(P)\bigcap W^{s}(Q)}$.  

Around each critical point in $I(P,R)$ we assume fixed a Morse chart inside which
the metric is the canonical one. Let $K_{\tau}(P,R)$
be the set of all points $x\in K(P,R)$ such that if
for some $t\in {\bf R}$ we have $\gamma_{t}(x)\in f^{-1}(c)$, then 
$d(\gamma_{t}(x),Q)\geq\tau$ for all $Q\in I(P,R)$ ($d( , )$ being the distance in
$M$). For $\tau$ sufficiently small this set is the (unreduced) suspension
over $Z_{\tau}(P,R)=(K(P,R)\bigcap f^{-1}(c))-\bigcup_{Q\in I(P,R)}D(Q,\tau)$ 
where $D(Q,\tau)$ is the disk in $M$ of center $Q$ and of radius
$\tau$  (the intersection of $S(Q,\tau)=\partial D(Q,\tau)$ and $f^{-1}(c)$ is 
certainly transverse for $\tau$ small enough).

Notice that $Z_{\tau}(P,R)$ is a manifold with boundary whose
homeomorphism type does not depend on the choice of $\tau$, if this constant is
smaller than some fixed $\tau'>0$, and that its interior is homeomorphic to $K(P,R)$. 
These statements follow from the Morse-Smale condition. Indeed, 
for $Q\in I(P,R)$ let $G(Q,\tau)=f^{-1}(c)\bigcap S(Q,\tau)$. Clearly,
$G(Q,\tau)\approx S^{q-1}\times S^{n-q-1}$. Consider a neighborhood of $Q$ as
described at the beginning of the proof of 3.1, $U'=U'_{\tau,\epsilon}$ and recall that
for small $\epsilon$ and $\tau$, 
$C'=\partial U'\bigcap f^{-1}(c+\epsilon)=D^{q}\times S^{n-q-1}$ and with this
identification $\{0\}\times S^{n-q-1}=S^{s}(Q)$. The intersection of $W^{u}(P)$
with $f^{-1}(c+\epsilon)$ is identified with $S^{u}(P)$ and therefore it intersects
transversely $S^{s}(Q)$ as well as $\partial C'$. This implies immediately that
$W^{u}(P)$ intersects transversely $G(Q,\tau)$. By the same method we
obtain that $W^{s}(R)$ intersects transversely $G(Q,\tau)$. This shows that for all
small enough $\tau$ the intersection of $K(P,Q)$ and $G(Q,\tau)$ is transverse. This
implies all the claimed properties of $Z_{\tau}(P,R)$.

Fix $Q$ and some $\tau$ as above and let $\partial_{Q} Z = \partial
Z_{\tau}(P,R)\bigcap S(Q,\tau)$, $\partial_{Q} K=\{x\in K_{\tau}(P,R): \exists t,
\gamma_{t}(x)\in S(Q,\tau)\bigcap f^{-1}(c)\}\bigcup \{P,R\}$. Clearly,
$\partial_{Q} K$ is the suspension over $\partial_{Q}Z$ and, in particular, if for $a\in
{\bf R}$ such that
$f(R)<a<f(P)$ we denote $\partial^{a}_{Q}=\partial_{Q}K\bigcap f^{-1}(a)$, then
$\partial^{a}_{Q}$ is homeomorphic to $\partial_{Q}Z$. 

\subsubsection*{Embedding of $\partial_{Q}K$ in $M$.}  Let
$W(P,Q)=\{x\in W^{u}(P) : f(x)\geq c, \exists t\ \gamma_{t}(x)\in f^{-1}(c),
d(\gamma_{t}(x),Q)\leq\tau\}\bigcup K(P,Q)$ and $W(R,Q) =\{x\in
W^{s}(R): f(x)\leq c, \exists t\ \gamma_{t}(x)\in f^{-1}(c),
d(\gamma_{t}(x),Q)\leq\tau\}\bigcup K(Q,R)$. 
Let $a$ and $b$ be such that $f(R)<b<c< a<f(P)$.
Clearly, the interiors of $W(P,Q)$ and $W(R,Q)$
are (open) cones over their intersection with $f^{-1}(a)$ and respectively
$f^{-1}(b)$. In
particular,  we have $W_{a}(P)=W(P,Q)\bigcap f^{-1}(a)\approx
D^{q}\times Z(P,Q)$ and $W_{b}(R)=
W(R,Q)\bigcap f^{-1}(b)\approx Z(Q,R)\times D^{n-q}$. The set
$H(Q,\tau)=\{x\in f^{-1}(c): d(x,Q)\leq\tau\}$ is clearly a cone over its
boundary which is $G(Q,\tau)$. Also $\partial^{a}_{Q}\subset W_{a}(P)$
and $\partial^{b}_{Q}\subset W_{b}(R)$. 

Notice, that the union $W(P,Q)\bigcup
W(R,Q)\bigcup H(Q,\tau)$ has the homotopy type of the wedge
$\Sigma Z(P,Q)\vee\Sigma Z(Q,R)$ and
$H(Q,\tau)\bigcup \partial_{Q} K\simeq \Sigma
\partial_{Q}Z\vee\Sigma\partial_{Q} Z$. Therefore the inclusion
$\partial_{Q} K\longrightarrow H(Q,\tau)\bigcup \partial_{Q} K \longrightarrow \\ 
  W(P,Q)\bigcup W(R,Q)\bigcup H(Q,\tau)\longrightarrow M$
is seen to be homotopic to 
$$
\label{eq:fact}
\partial_{Q}K=\Sigma\partial_{Q}Z\longrightarrow
\Sigma\partial_{Q}Z\vee\Sigma\partial_{Q}Z\longrightarrow \Sigma Z(P,Q)\vee
\Sigma Z(Q,R)\longrightarrow M 
$$
where the first map is the pinch map, the second is given by wedging the
suspensions of
$t_{1}:\partial^{a}_{Q}\subset W_{a}(P)=D^{q}\times
Z(P,Q)\stackrel{p_{2}}{\longrightarrow}Z(P,Q)$ and
$t_{2}:\partial^{a'}_{Q}\subset W_{a'}(R)=Z(Q,R)\times
D^{n-q}\stackrel{p_{1}}{\longrightarrow}Z(Q,R)$ and the third map is the
inclusion of $K(P,Q)\vee K(Q,R)$.

Because $M$ is simply connected we have a well defined map 
$l(P,R):\partial_{Q}Z\longrightarrow\Omega M$ which is the
adjunct of the inclusion $\partial_{Q}K\subset M$. From the factorization above,
by adjunction, we obtain that $l(P,R)$ factors as 
\begin{equation}
\label{eq:fact}
\partial_{Q}Z\longrightarrow Z(P,Q)\times Z(Q,R)\longrightarrow \Omega M\times
\Omega M\longrightarrow \Omega M
\end{equation}
where the first map is the product $t_{1}\times t_{2}$, the second map is the
product of the adjoints of the respective inclusions and the third is the loop product. 

\subsubsection*{Description of $\partial_{Q}Z$.} We want to observe that
$t_{1}\times t_{2}$ is a homeomorphism. Fix $a$ and $b$ such that $f(R)<b\leq
c\leq a<f(P)$ and $c-b<\epsilon$, $a-c<\epsilon$. We have the
homeomorphsims
$\partial^{b}_{Q}\approx\partial_{Q} Z\approx\partial^{a}_{Q}$ induced by
transporting the respective sets along the flow
$\gamma$. Consider the inclusion $\partial^{a}_{Q}\hookrightarrow S^{q-1}\times
Z(P,Q)$ (here $S^{q-1}\times Z(P,Q)\subset D^{q}\times Z(P,Q)= W^{u}(P)\bigcap
f^{-1}(a)\bigcap U'$). By transporting it along the flow $\gamma$ till we reach
the level set $f^{-1}(b)$ we obtain an inclusion $\partial^{b}_{Q}
\hookrightarrow S^{q-1}\times Z(P,Q)\stackrel{j}{\hookrightarrow} S^{q-1}\times
S^{n-q-1}=\partial U'\bigcap f^{-1}(b)$ (notice that now $S^{q-1}$ is identified here
to $S^{u}(Q)$). The map $j$ is the identity on the first factor and the usual inclusion
on the second. At the same time we also have an inclusion $\partial^{b}_{Q}
\hookrightarrow Z(Q,R)\times S^{n-q-1}$. It follows that the product $t_{2}\times
t_{1}$ is identified to the inclusion $v:\partial^{b}_{Q}\hookrightarrow
Z(Q,R)\times Z(P,Q)=(Z(Q,R)\times S^{n-q-1})\bigcap (S^{q-1}\times Z(P,Q))$.
On the other hand each point of this last intersection belongs to $\partial^{b}_{Q}$.
Indeed, with the given parametrizations all the points in $S^{q-1}\times Z(P,Q)$
belong to $W^{u}(P)$ and all the points in $Z(Q,R)\times S^{n-q-1}$
belong to $W^{s}(R)$. Therefore $v$ is also surjective.

\subsubsection*{Consequences in $\Omega_{\ast}^{fr}(\Omega M)$.}
Consider the framing of $\partial_{Q}Z=Z(P,Q)\times Z(Q,R)$ inside
$S^{n-q-1}\times S^{q-1}$ obtained as the restriction of the normal framing
of $K(P,R)$ inside $f^{-1}(a)$. This framing is the product of the
framing induced by the standard framing of $S^{p-1}$ with
that induced by the standard framing of $S^{n-r-1}$. Therefore, the framing
of $Z(P,Q)\times Z(Q,R)$ inside $S^{n-q-1}\times S^{q-1}$ coincides with
the product of the standard framings of $Z(P,Q)$ in $S^{n-q-1}$
(induced by the framing of $S^{p-1}$) and 
the standard framing of $Z(Q,R)$ in $S^{q-1}$. 
Let $\overline{Z}(P,Q)$ be the framed bordism
representative given by $l(P,Q)$ and the framing of $Z(P,Q)$ inside $S^{n-q-1}$.
The orientation of $\overline{Z}(P,Q)$ is such that the framing
summed with this orientation gives the standard orientation on $S^{n-q-1}$.

We orient $Z(P,R)$ such that the
ordered sum of the framings induced first from $S^{p-1}$, from
$S^{n-r-1}$ and the orientation of $Z(P,R)$ gives the orientation of $f^{-1}(a)$.

Assume now that $M$ has trivial stable normal bundle. In this case, 
after embedding $M$ in a high dimensional sphere we see that
because of (\ref{eq:fact}) and as
 $$\coprod_{Q\in I(P,R)}\partial_{Q}Z=\partial
(Z_{\tau}(P,R))$$ we have $\Sigma_{Q\in I(P,R)} \epsilon
'(Q) [\overline{Z}(P,Q)]^{fr}\bullet [Z(Q,R)]^{fr}=0$. 

The sign $\epsilon '(Q)$ is $+1$
if the orientation induced from that of $Z(P,R)$ on $\overline {Z}(P,Q)\times
Z(Q,R)$ coincides with the product orientation and is $-1$ otherwise.
With our conventions $\epsilon'(Q)=(-1)^{(n-r)q}\epsilon_{0}$
where $\epsilon_{0}$ depends only on $p,r,n$.

We now return to the case when the stable bundle of $M$ is general.
Similarly to the proof of Proposition 3.2 we consider a function
$f'=f\circ p+g$ defined on a neighborhood $U$ of the $0$-section of the normal
bundle $\nu$  of an embedding of $M$ in a high dimensional sphere $S^{n+m}$.
Here, $g$ is the square of a distance function measuring the distance from $M$. 
We use the notation  $Z^{-\nu}(P,Q)$ for the connecting manifold of $P$ and $Q$
for  the function $f'$. By the formula above we have 
$\Sigma_{Q\in I(P,R)}(-1)^{(m+n-r)q}[\overline{Z^{-\nu}}(P,Q)]^{fr}\bullet
[Z^{-\nu}(Q,R)]^{fr}=0$. We have $[Z^{-\nu}(Q,R)]^{fr}=[Z(Q,R)]^{fr}$,
$[Z^{-\nu}(P,Q)]^{fr}= [Z(P,Q)]^{fr}$.  
Moreover, because the stable normal bundle of $U$ is trivial, 
the framings of $\overline{Z^{-\nu}}(P,Q)$ and
$Z^{-\nu}(P,Q)$ are the same up to a sign coming from a possible difference
in orientations. With our conventions the sign is $(-1)^{(m+n-p)q}$.
\end{proof}

\

\begin{rem} In \cite{CohenJS},\cite{CohenJS2} Cohen, Jones and Segal pursue a
systematic analysis of the ends of the moduli spaces of connecting flow lines
(even for non consecutive critical points) together with the relevant framings (the
definition they use for these framings is somewhat different from ours, though). 
In particular,  the fact that
$t_{1}\times t_{2}$ above is a homeomorphism is a direct consequence of a more
general result of Betz as described in \cite{CohenJS}.  For completeness we have
included a direct justification in the proof above. 
\end{rem}

\section{Applications and Examples.}

As in the last section $M^{n}$ is a compact, smooth, riemannian,
simply connected manifold. Suppose $f:M\longrightarrow {\bf R}$ is a Morse-Smale
function and
$P$ and $Q$ are consecutive critical points of index, respectively,
$p$ and $q$. Recall that
$h^{\ast}(P,Q)\in\pi_{p-q-1}\Omega^{q}\Sigma^{q}\Omega M$ is the $q$-th
order adjoint of $h(P,Q)$. Let
$h^{S}(P,Q)\in\pi_{p-q-1}^{S}(\Omega M)$ be the class of
$h^{\ast}(P,Q)$. Let $k(P,Q)=min\{ j : h^{S}(P,Q)\in Im (\Omega^{j}\Sigma^{j}
\Omega M)\}$. This gives a measure of the embedding complexity of $Z(P,Q)$.
Indeed, if the framed embedding $Z(P,Q)\hookrightarrow S^{p-1}=S^{u}(P)$
factors as a composition of framed embeddings $Z(P,Q)\hookrightarrow
S^{t-1}\hookrightarrow S^{p-1}$, then $k(P,Q)\leq t-p+q$. For example, if
$k(P,Q)>0$, then $Z(P,Q)$ is not a sphere. 

Of course, as $k(P,Q)$ is defined homotopically, it is invariant to deformations
of $f$ via Morse-Smale functions. 

We intend here to give a method to estimate $k(P,Q)$.
In particular, we construct examples when $k(P,Q)=ind (Q)$ even if
$\delta_{f}(P,Q)=0$.  Of course, we always have $k(P,Q)\leq ind(Q)$.
Along the way, we also detect certain
Poincar\'e complexes that are not smoothable.

\subsection{Non-smoothable Poincar\'e duality complexes and embedding complexity
of connecting manifolds.}

We start with two simple consequences of 3.2. 
\begin{cor}
Assume that $f:M\longrightarrow {\bf R}$ is a smooth
Morse-Smale function and that $P$ and $Q$ are consecutive critical points of
$f$. Then, for some $\epsilon\in\{-1,1\}$, we have $\epsilon \{\delta_{-f}(Q,P)\}
=\{\delta_{f}(P,Q)\}+\Delta (P,Q)$. Where $\{x\}$ is the stable image of $x$ and
$\Delta (P,Q)=J^{q}(\Omega^{q}\Sigma^{q}(\Omega \nu)\circ h^{\ast}(P,Q))$
($\nu$ is the stable normal bundle of $M$).
\end{cor}

\begin{proof} We apply 3.2 to the stable normal bundle of $M$ together with
the duality results (mentioned in 2.3) implying $\delta_{-f}(Q,P)=\epsilon
\delta_{f}(P,Q)^{\nu}$. 
\end{proof}

Suppose now that $H_{\ast}(M;{\bf Z})$ is torsion free and that for some $q<p<n/2$
we have $H_{q}(M;{\bf Z})\approx H_{p}(M;{\bf Z})\approx {\bf Z}$,
$H_{\ast}(M;{\bf Z})=0$ for $q<\ast<p, \ast=q-1, \ast=p+1$. 
Then, in a minimal cell
decomposition of $M$, there are two pairs of dual cells $e^{p}$, $e^{n-p}$
and $e^{q}$ and $e^{n-q}$ representing Poincar\'e dual generators in
$H_{p}(M)$, $H_{n-p}(M)$ and, respectively, $H_{q}(M)$ , $H_{n-q}(M)$.
As $e^{p}$, $e^{q}$ are attached in succession there is a relative
attaching map $\delta : S^{p-1}\longrightarrow S^{q}$ and similarly
$e^{n-q}$ and $e^{n-p}$ being attached successively we have another relative
attaching map $\delta' : S^{n-q-1}\longrightarrow S^{n-p}$.

\

In all this section $\epsilon\in \{-1,+1\}$. Its presence in the formula below
reflects an indeterminancy caused by the fact that
if $\delta$ is the relative attaching map above, then there is a different minimal
cell decomposition having as corresponding relative attaching map $-\delta$.

\begin{cor} In the setting above $\delta$ and $\delta'$ are (up to sign) independent
of the minimal $CW$ -decomposition used in their definition and
$\{\delta\}=\epsilon\{\delta'\} $ mod  $ (Im(J^{q}))$.
\end{cor}

\begin{proof} The first part is immediate.We apply the corollary above to a perfect
Morse function and its consecutive critical points $P$, $Q$ of indexes 
$p$ and $q$. Because of the first part we have $\delta_{f}(P,Q)=\delta$ and
$\delta_{-f}(Q,P)=\delta'$.
\end{proof}

\begin{rem}  a. Corollary 4.2 can also be proven by purely homotopical methods.
Surprisingly, it appears not to have been known before.

b. For spaces more general than those appearing in 4.2
the relative attaching maps of the type of $\delta$ and $\delta'$ 
depend on the specific minimal cell-decomposition used. 
Here is a relevant example
(well known in the study of the Mislin genus \cite{McGibbon}). 
Let $X=S^{9}\bigcup_{5\nu}e^{13}$ and
$X'=S^{9}\bigcup_{13\nu}e^{13}$ here $\nu\in \pi_{12}(S^{9})\approx \pi^{S}_{3}
\approx {\bf Z}/24$ is a generator. These two complexes are not homotopy equivalent
(even if they are so when localized at each prime $p$). However, it is easy  to show
that $X'\vee S^{9}\simeq X\vee S^{9}$. Hence, we have two very different minimal 
cell decompositions for the space $Y=X\vee S^{9}$. In particular, if $N$ is a
manifold with simply connected boundary having the homotopy type of $Y$ and
of dimension greater than $26$, each of these two cell decompositions is induced by a
perfect Morse-Smale function. These two Morse-Smale functions are in different
connected components of the space of perfect Morse-Smale functions of $N$.
However, they are in the same component of the space of perfect {\em Morse}
functions. Indeed, by a result of Matsumoto \cite{Matsumoto} the space of perfect
Morse functions on a simply connected manifold with simply connected boundary, of
dimension greater than $5$ and with torsion free homology is connected
(see also \cite{Agoston}, for more general results \cite{Sharko}).
\end{rem}

\

It is immediate to see that if $\Delta (P,Q)
\not\in Im(J^{t}) $, then $k(P,Q)>t$. The obvious question that we
consider now is what values can take the twisting $\Delta (P,Q)$. 

\begin{rem}
Certainly, there are examples when $\Sigma H(P,Q)$ is not vanishing, but
$\Delta (P,Q)$ is null. An instructive example
is that of ${\bf CP}^{n}$ together with a perfect Morse-Smale function.
If $P$ and $Q$ are two arbitrary consecutive critical points of $f$, then
$p=ind(P)=ind(Q)+2$ and $H(P,Q)\in\pi_{p-1}(\Sigma^{p-2}\Omega {\bf
CP}^{n})\approx {\bf Z}$ is a generator. On the other hand the
difference $\delta_{f}(P,Q)-\delta_{-f}(Q,P)$ is non zero iff $n$ is even
(because in that case $w_{2}\not= 0$). 
\end{rem}

\

In the constructions below the strategy will be the following:
\begin{itemize}
\item Construct a certain Poincar\'e duality $CW$-complex $X$ with an explicit cell
decomposition.
\item Assuming $X$ smooth, consider a perfect Morse-Smale function $f$ on $X$
which induces the fixed cell-decomposition at least below the middle dimension.
\item Evaluate the restriction of the stable normal bundle of $X$  to the middle
dimensional skeleton.
\item Identify two consecutive critical points $P$ and $Q$ and evaluate
using 4.1, 4.2 or 3.2 the twisting $\Delta (P,Q)$.
\end{itemize}

\subsubsection*{Some examples of Larry Smith.} 
Let $X=(S^{p}\vee S^{n-p})\bigcup_{h} e^{n}$ where $h=[i_{1},i_{2}]+i_{2}\circ
x$ with $p<n/2$, $i_{1}:S^{p}\longrightarrow S^{p}\vee S^{n-p}$ and
$i_{2}:S^{n-p}\longrightarrow S^{p}\vee S^{n-p}$ the inclusions, $[-,-]$ the
Whitehead product, $x\in\pi_{n-1}(S^{n-p})$ such that $x\not\in Im (J)$.
Then, obviously $X$ is a Poincar\'e duality space, but it is not smoothable
by 4.2 as the stable difference of the two relative attaching maps associated
the first to the $p$ and $0$ cells, and the second to the $n$ and $n-p$ cells, is 
equal to $x$.

Smith has proved this result by purely homotopical methods \cite{Smith}. We
indicate below a different, very short, purely homotopical proof of this same fact. It
pinpoints the homotopical content of our Morse theoretical techniques. 

Assume that $X$ is smoothable. Thus, its stable normal bundle $\nu$ is orthogonal.
Fix $m=rank (\nu)$ and let $T^{\nu}(X)$ be the associated Thom space. It is
Spanier-Whitehead dual to $X_{+}$. In particular, the first notrivial relative attaching
map in $T^{\nu}(X)$, $\delta:S^{m+p-1}\longrightarrow S^{m}$ equals stably
(up to sign) the top relative attaching map of $X$, $\delta':S^{n-1}\longrightarrow
S^{n-p}$. Of course, $\delta'=x$ hence $\delta=\epsilon x$. Now, $\delta$ is already
present in the Thom space, $T^{\nu}(S^{p})$, of the restriction of $\nu$ to the bottom
sphere $S^{p}$. By classical results (appearing already in \cite{Milnor}
or \cite{Adams}) we have that $T^{\nu}(S^{p})\simeq S^{m}\bigcup _{\tau}
e^{m+p}$ and $\tau\in Im(J)$. As $\tau=\delta$ we are led to a contradiction.

\begin{rem} 
a. Of course, if $x\in Im(J_{PL})$, then $X$ has the homotopy type of a $PL$
manifold \cite{Smith}.  

b. Clearly, one may extend the examples
above by producing  $k-1$- connected Poincar\'e duality
complexes of dimension $n$, $2\leq k<n/2$ with a relative
attaching map of the top cell with respect to some $n-k$ cell
that does not belong to $Im(J)$. Of course, this is the content of Spivak's
"first smoothing obstruction" \cite{Spivak}( see also \cite{KervaireM}) which is thus
recovered from 4.1. Notice, on the other hand, that Kervaire's original
non-smoothable PL-manifold \cite{Kervaire}, as well as the examples of
Eells and Kuiper \cite{EellsK}, are not detected by these means. Indeed, the
non-smoothable Poincar\'e spaces obtained by our methods are all stable in the
sense that they remain non-smoothable after crossing with a sphere.
\end{rem}

\subsubsection*{Non-smoothable Poincar\'e spaces with vanishing Spivak first
smoothing obstruction.} We construct here an example of a non-smoothable
Poincar\'e complex whose non-smoothability can not be detected by 
the relative attaching maps of the top cell.

Take $q>2$ and
$Z=(S^{2}\vee (S^{q}\bigcup_{\eta}e^{q+2}))\bigcup_{[i,j]}e^{q+4}$ with
$i:S^{2}\longrightarrow S^{2}\vee (S^{q}\bigcup_{\eta}e^{q+2})$ the inclusion,
$\eta\in\pi_{q+1}(S^{q})\approx {\bf Z}/2$ a generator and
$j\in \pi_{q+2}(S^{q}\bigcup_{\eta}e^{q+2})$ such that the image of
$j$ via the pinching map $p: S^{q}\bigcup_{\eta}e^{q+2}\longrightarrow S^{q+2}$
is equal to twice a generator of $\pi_{q+2}(S^{q+2})$. Consider
$y\in\pi^{S}_{q+1}$ such that $2y\not\in Im(J^{2})$. Let $BF$ be the classifying
space of stable spherical fibrations. As $\pi_{k}(BF)=\pi_{k-1}^{S}$
there is a spherical fibration given by the composition
$\mu: Z\stackrel{t}{\longrightarrow}S^{q+2}\stackrel{u}{\longrightarrow}
BF_{m}\stackrel{v}{\hookrightarrow} BF$ with $y\simeq (v\circ u)^{\ast}$,
$BF_{m}$ is the classifying space of spherical fibrations of fibre $S^{m-1}$ ($m$ big
enough) and $t$ is induced by the collapsing $S^{2}\vee
(S^{q}\bigcup_{\eta}e^{q+2})
\longrightarrow S^{q+2}$.  Let $N$ be a trivial, smooth thickening of $Z$  - this
is a smooth manifold with boundary of dimension $n>2(q+4)+1$
having the homotopy type of $Z$ and which embedds in $S^{n}$ 
 (by \cite{Wall} it is unique up to isotopy). Let $N'$ be the double of $N$.
As $N'$ is the union of two copies of $N$ pasted over $\partial N$, two copies of the
spherical fibration $\mu$ can be also pasted together to give a new spherical fibration
$\mu':N'\longrightarrow BF_{m}$. The total space $K$ of $\mu'$ is our example.

\begin{rem} There is another more elegant way to produce $\mu'$ out of $\mu$.
It was indicated to me by Mark Mahowald. 

First a general remark of some independent interest (it appears to be
folklore, for extensions see \cite{Klein2}). Because
$N'$ is the double of a trivial thickening, the suspension $\Sigma N'$ splits
as $\Sigma N\vee \Sigma^{2} N^{\#}\vee S^{n+1}$ with $N^{\#}$ homotopy
equivalent to the complement of $N$ in $S^{n}$. Indeed, $N$ is embeddable
in $S^{n}$, therefore we have $S^{n}=N\bigcup_{\partial N}N^{\#}$. The 
inclusions of $N$ and $N^{\#}$ in $S^{n}$ being homotopically trivial, we obtain a
cofibration sequence $S^{n}\longrightarrow \Sigma\partial N\longrightarrow \Sigma
N\vee\Sigma N^{\#}$. This cofibration sequence splits and thus
$\Sigma\partial N \simeq \Sigma N\vee\Sigma N^{\#}\vee S^{n}$. Now
$N'=N\bigcup_{\partial N}N$ hence $\Sigma N'=\Sigma N\bigcup_{\Sigma\partial
N}\Sigma N$. The description of $\Sigma\partial N$ implies the claimed splitting.

Consider $\mu^{\ast}:\Sigma N\longrightarrow B^{2}F$ whose adjoint is $\mu$.
It extends trivially to $(\mu^{\ast})':\Sigma N'=\Sigma N\vee\Sigma^{2}N^{\#}\vee
S^{n+1}\longrightarrow B^{2}F$. Its adjoint gives $\mu'$. 
\end{rem}

\begin{lem}
The space $K$ does not have the homotopy type of a smooth manifold.
\end{lem}

\begin{proof}
Assume $K$ smooth. Its stable normal
bundle $\nu:K\longrightarrow {\bf BSO}(m')$, $m'>m+n$, has the property that
if $\nu':K\longrightarrow {\bf BSO}(m')\longrightarrow BF_{m'}$ is
the associated spherical fibration, then $\nu'|_{N}+\mu=0$ in $[N, BF]$ (because
this sum gives the Spivak spherical fibration of $N$ which is trivial). Therefore
$\nu'|_{N}=-\mu$. Let $f:K\longrightarrow {\bf R}$ be a perfect Morse-Smale
function inducing the standard cell decomposition given for $Z$ on the $q+4$
dimensional skeleton of $K$, $K^{(q+4)}\simeq Z$. We may assume also that $f$ is
such that the critical points $P$ and $Q$ of indexes, repectively, $q+4$ and $2$ are
consecutive. Of course, $\delta_{f}(P,Q)=0$. We intend to evaluate
$\Delta (P,Q)$. Consider the following diagram.
$$ \xy
\xymatrix@+2pt{
S^{q+3}\ar[r]\ar[dr]_{l}\ar[ddr]_{s}&S^{2}\wedge \Omega (S^{2}\vee
(S^{q}\bigcup e^{q+2}))\ar[d]^{ S^{2}\wedge\Omega\nu}\ar[r]&
S^{2}\wedge\Omega S^{q+2}\ar[ddl]^{-\Sigma^{2}\Omega u}\\
&S^{2}\wedge {\bf SO}(m')\ar[d]&\\
&S^{2}\wedge \Omega BF_{m'}&}
\endxy $$

By the basic properties of the Whitehead product, the top horizontal composition
is twice the inclusion of the bottom cell. This implies that $s=-2\Sigma^{2}y$.
Using 4.1 and the diagram below we see that $\delta_{-f}(Q,P)$ 
agrees stably (up to sign) with the composition
$k\circ \Sigma^{m'}l=k'\circ\Sigma^{m'}s$ .
$$ \xy
\xymatrix@+1cm{
S^{m'}\wedge S^{q+3}\ar[r]^{\Sigma^{m'}l}\ar[dr]_{\Sigma^{m'}s}&
S^{m'}\wedge S^{2}\wedge{\bf SO}(m')\ar[r]^{\ \ k}\ar[d]&S^{m'+2}\ar[d]^{id}\\
& S^{m'}\wedge S^{2}\wedge\Omega BF_{m'}\ar[r]_{\ \ k'}&S^{m'+2}}
\endxy $$

Here
$k$ and $k'$ are the double suspensions of the maps induced by the Hopf
construction applied respectively to the actions $S^{m'-1}\times{\bf
SO}(m')\longrightarrow S^{m'-1}$ and $S^{m'-1}\times
\Omega BF_{m'}\longrightarrow S^{m'-1}$. As $s$ is itself a double suspension we obtain
that $k'\circ \Sigma^{m'}s$ is stably equal to the class of $s$ in $\pi_{q+1}(\Omega
BF)=\pi_{q+1}^{S}$. Therefore, $\Delta (P,Q)= -2y$ which  leads to a contradiction.
\end{proof}

Clearly the relative attaching maps of the top cell $S^{n+m-1}\longrightarrow
S^{n+m-2}$ and $S^{n+m-1}\longrightarrow S^{n+m-q}$ vanish.

\begin{rem} a. It is easy to see that if $y$ lifts to an element in
$\pi_{\ast}(\Omega BPL)$, then $\mu$ and $\mu'$  are PL-bundles and therefore $K$
is a PL-manifold. 

b. Mark Mahowald has described another interesting non-smoothable Poincar\'e
space that we now present (see also \cite{MahowaldT} p.408). 

Let $U$ be the total space of the spherical fibration classified by
$\eta^{\ast}_{4}\in\pi_{17}(BF_{14})$
 where $\eta^{\ast}_{4}$ is the adjoint of the fourth element of
Mahowald's family  \cite{Mahowald2}. Recall that $\eta_{4}\not\in Im(J)$.

A cell decomposition of this space is as follows.
Let $i_{n}$ be a generator of $\pi_{n}(S^{n})$ and $\nu_{n}$ a generator
of $\pi_{n+3}(S^{n})$. Let
$T=S^{13}\bigcup_{\nu_{13}}e^{17}$. The Whitehead product
$[\nu_{13},i_{13}]$ being null (\cite{Toda}, \cite{Mahowald1}) the
relative Whitehead product $[i_{17},i_{13}]\in\pi_{29}(T, S^{13})$ pulls back to an
element $z\in\pi_{29}(T)$ and $U=T\bigcup_{z}e^{30}$. 

Notice that $\Sigma^{n}U$ splits as
$S^{13+n}\bigcup_{\nu_{13+n}}e^{17+n}
\bigcup_{i_{13+n}\circ\eta_{4}}e^{n+30}$.
 
The space $U$ is not smoothable. 
This is detectable by a top cell argument as follows.

We have $T=\Sigma T'$. Thus, we may interpret
$T$ as the cofibre of a trivial attaching map $T'\longrightarrow\ast$.
Therefore, $U\times S^{n}$ has a cone decomposition that starts with $T\vee S^{n}$
and has two other stages given by cofibration sequences $T'\ast
S^{n-1}\longrightarrow T\vee S^{n}\longrightarrow U'$ and $S^{29}\ast
S^{n-1}\longrightarrow U' \longrightarrow U\times S^{n}$. 
By a Spanier-Whitehead duality argument, one sees that if $U$ and hence $U\times
S^{n}$ are smoothable, then the relative attaching map $S^{n+29}\longrightarrow
\Sigma^{n}T$ should not contain $\eta_{4}$. 

Obviously,  this example does not fit our Morse theoretic setting.
However, a critical point approach is possible by using degenerate critical
points as in \cite{Cornea2}.
\end{rem}

\subsubsection*{Morse-Smale functions with twisting realizing elements of $Im
(J^{q})$.} We consider here the question of what elements in $Im(J^{q})$ can be
realized as differences, $\Delta (P,Q)$, of attaching maps like in Corollary 4.1.

\begin{lem} For any $x\in Im(J^{q})$, $q\geq 2$, there is a 
smooth manifold $M$ and a Morse-Smale function $f:M\longrightarrow {\bf R}$
with consecutive critical points $P$ and $Q$ such that $Q$ has index $q$
and $\delta_{f}(P,Q)=0$, $\delta_{-f}(Q,P)=\epsilon x$.
\end{lem}

\begin{proof} As $x\in Im(J^{q})$ there is $x^{1}:S^{p-1}\longrightarrow
\Sigma^{q} {\bf SO}(m)$ such that $x$ is given by $J^{q}(x^{\ast})$ where
$x^{\ast}$ is the $q$-th order adjoint of $x^{1}$ and $m$ is large enough. 
Denote by $BS'$ the $p-q$-dimensional skeleton of ${\bf BSO}(m)$ and let
$S'=\Omega BS'$. Of course, $x^{1}$ lifts to an element $x' \in
\pi_{p-1}(\Sigma^{q}S')$.  

There is a
fibration $\Omega S^{q}\ast S'\stackrel{j}{\longrightarrow} S^{q}\vee
BS'\longrightarrow S^{q}\times BS'$ and a natural
inclusion $\Sigma^{q} S'=S^{q-1}\ast S'\stackrel{i}{\hookrightarrow} \Omega
S^{q}\ast S'$. Let $x''=j\circ i\circ x'$. Denote by $L$ the cofibre of $x''$.
Fix a smooth, trivial, $n$-dimensional thickening of $L$ denoted by $H$ and let $H'$
be its (smooth) double. We need to define a certain fiber bundle over $H'$. Denote by
$\mu:H\simeq L\longrightarrow {\bf BSO}$ the bundle defined by the 
trivial extension
of the map $S^{q}\vee BS' \hookrightarrow
S^{q}\vee{\bf BSO}(m)\stackrel{c}{\longrightarrow} {\bf BSO}(m)\hookrightarrow
{\bf BSO}$ ($c$ being the obvious collapsing). As discussed above one can extend
$\mu$ to a bundle $\mu':H'\longrightarrow {\bf BSO}$ and it lifts to a bundle 
 $\tau :H'\longrightarrow {\bf BSO}(m')$ with $m'>n,m$.
Let $G$ be the total space of the associated sphere bundle. As $H'$ is smooth
we see, by a transversality argument, that we may assume $\tau$ smooth
and therefore $G$  is also smooth.  There is a perfect Morse-Smale function
$f:G\longrightarrow {\bf R}$ such that the critical point $P$ of index
 $p$ and the critical point $Q$ of index $q$ corresponding to $S^{q}\hookrightarrow
L\hookrightarrow H'$ are consecutive, and such that the attaching map of
the $p$-dimensional cell is given by $x''$. We have $\delta_{f}(P,Q)=0$.
We now intend to compute $\delta_{-f}(Q,P)$.  Let $\nu$ be the stable normal
bundle of $G$. Obviously, its restriction to $L\subset H'$ coincides with $-\mu$.
Notice also that, as $\mu$ is null on $S^{q}$, 
$\Delta (P,Q)$ can be estimated by using the composition
$S^{p-1}\stackrel{x'}{\longrightarrow}S^{q}\wedge S'\hookrightarrow
S^{q} \wedge\Omega G$ instead of $H(P,Q)$. Indeed, the image of this composition
by the map $S^{q}\wedge \Omega G\longrightarrow S^{q}\wedge {\bf SO}(m')$
coincides with the image of $H(P,Q)$.
The composition $S^{q}\wedge S' \hookrightarrow
 S^{q}\wedge\Omega G\stackrel{\Sigma^{q}\nu}{\longrightarrow} S^{q}\wedge{\bf
SO}(m')$ is homotopic to the negative of the inclusion $ S^{q}\wedge S'
\hookrightarrow  S^{q}\wedge {\bf SO}(m')$.

Therefore, by applying 3.2, we see that $\delta_{-f}(Q,P)$ is (up to sign)
equal to $J^{q}(x^{\ast})=x$. 
\end{proof}

\begin{rem} It is likely that the Corollaries 4.1 and 4.2 as well as Proposition 3.2
 have analogues in the $PL$ and $Top$ categories. 
\end{rem}

\subsection{An extension of the Morse complex and detection of connecting
flow lines.}

The results of the second section can be used for
the detection of critical points and connecting flow lines. Indeed, assume
that $P$, $Q$ are consecutive critical points of the
Morse-Smale function $f:M\longrightarrow {\bf R}$. If the
suspension of the Hopf invariant $H(P,Q)$
is not zero, then there are flow lines connecting $P$ to $Q$. Of course, such
connecting flow lines are already known to exist if $\delta_{f}(P,Q)\not=0$
but, as seen before, this relative attaching map can vanish without
$\Sigma H(P,Q)$ being also null.

\subsubsection*{The energy functional on $\Omega S^{3}$.} This example has
been suggested to me by Raoul Bott. Consider the energy functional $E$ on 
$\Omega (S^{3};u,v)$,  the space of piecewise smooth curves connecting the points
$u, v\in S^{3}$,  (with $u$ and
$v$ different and not antipodals).  It is well known that $E$ is a perfect Morse function
\cite{Milnor3} whose critical points are the geodesics connecting $u$ and $v$. 
Fix two such geodesics  $P$  and $Q$. 
Approximate $\Omega (S^{3}; u,v)$ by a finite dimensional manifold that contains
$P$ and $Q$, on which the restriction of $E$, $E'$, is still Morse, has
critical points of the same index and the sets $\Omega ^{a}=(E')^{-1}(-\infty, a]$ are 
compact for $a\leq b=\max\{E(P),E(Q)\}+\tau$ ($\tau >0$ small) and
have the same homotopy type as $E^{-1}(-\infty, a]$ \cite{Milnor3}. Fix a metric
on $\Omega ^{b}$ such that $E'$ is Morse-Smale. We want to remark that $P$ and
$Q$ are connected  by some flow lines of the flow induced by $-\nabla E'$. 
Notice that $\Omega (S^{3}; u,v)$ is homotopy equivalent to $\Omega S^{3}$.
Assume $R$ is another critical point of $E$ in $\Omega^{b}$ such that
$ind(P)=ind(R)+2$, then
$\delta_{f}(P,R)=0$. However,
$\Sigma H(P,R)\not=0$. This implies, by transitivity, that $P$ and $Q$  are
also connected by some flow line. 

\subsubsection*{Extension of the Morse complex.}
The "detection" arguments above can be pursued further. For example, if 
$P$, $Q$ and $Q$, $R$ are two pairs of consecutive critical points 
for which the homology classes 
$[h(P,Q)], [h(Q,R)]\in
H_{\ast}(\Omega M;{\bf Z}/2)$ are known and if $[h(P,Q)]\bullet [h(Q,R)]\not= 0$
 it follows from 3.3 that there is at least another critical point $Q'$ that is connected via
possibly broken flow lines to both $P$ and $R$ ($\bullet$ being here the Pontryagin
product). 

Here is a way to encode in a somewhat global fashion the type of information
given by 3.3.

As before $f:M\longrightarrow {\bf R}$ is a Morse-Smale function, constant, regular
and maximal on $\partial M$. Let $a_{1}<a_{2}<...< a_{r}$ be a set $A$ of real
numbers such that the critical values of $f$ all appear among the
$a_{i}$'s as well as $f(\partial M)$. Let $\mathcal{C}_{k}=\Omega
_{\ast}^{fr}(\Omega M)<X\in f^{-1}(a_{k}) : \nabla f(X)=0>$ (where
$R<X_{1},...,X_{n}>$ is the free $R$-module generated by
$X_{1}$, ..., $X_{n}$). If $X$ is a critical point let $x$ be its index. Let
$d: \mathcal{C}_{i}\longrightarrow \mathcal{C}_{i-1}$ be the unique
$\Omega^{fr}_{\ast}(\Omega M)$-module morphism given by $$d(P)=\sum_{X\in
f^{-1}(a_{i-1}),\nabla f (X)=0}(-1)^{px}[Z(P,X)] X$$

\begin{cor}
The graded $\Omega^{fr}_{\ast}(\Omega M)$-module $(\mathcal{C}_{\ast}, d)$
is a chain complex. 
\end{cor}

\begin{proof} Notice that any pair of critical points $y\in f^{-1}(a_{i})$, 
$x\in f^{-1}(a_{i-1})$ are consecutive and apply 3.3.
\end{proof}

\begin{rem} a. 
A somewhat easier to handle complex is obtained by replacing
$\Omega_{\ast}^{fr}(\Omega M)$  with the Pontryagin ring
$H_{\ast}(\Omega M; {\bf Z}/2)$ and using instead of $d$,
$d'(P)=\sum_{X\in f^{-1}(a_{i-1}),\nabla f (X)=0}H'(P,X)X$
where, for two consecutive critical points $P$ and $Q$, we denote by
$H'(P,Q)$ the number of elements in $Z(P,Q)$ if $p=q+1$ and the
homology class of $H(P,Q)$ in $H_{p-q-1}(\Omega M; {\bf Z}/2)$ if $p>q+1$.
If $f$ is a self indexed function and $A={\bf N}$, then this complex
is the Morse complex of $f$ tensored with $H_{\ast}(\Omega M;{\bf Z}/2)$
(recall that $M$ is simply connected). 

b. It would be interesting to know whether one can deduce the existence of the
discussed chain complex by analytic methods.
\end{rem}

\subsubsection*{Fusion of critical points.}
One initial motivation for this work was the problem of
constructing functions with the least possible number of (possibly degenerate)
critical points on a given smooth manifold $M$. It is well-known
that a strict lower bound for this number is the Lusternik-Schnirelmann category
$cat(M)$ of $M$ \cite{LusternikS}. It was shown in \cite{Cornea4} that
when $M$ is $2$-connected, and $k\geq dim(M)$ there is a function on 
$M\times D^{k}$ regular, maximal and constant on $\partial (M\times D^{k})$
which realizes the lower bound given by the category up to one unit. The next
step is to understand when two consecutive critical points of a given function can be
"fused" together. More precisely, let $f:M\longrightarrow {\bf R}$ be a smooth
function and let $P$, $Q$ be consecutive critical points of $f$. The question
is whether there is a function $f'$ equal to $f$ in the exterior of some neighborhood
$U$ of the closure of the points situated on flow lines connecting $P$ to $Q$
and having just one critical point in $U$. It is natural to first assume that $P$
and $Q$ are consecutive, non-degenerate critical points of indexes, respectively, 
$p$ and $q$, and
that, with respect to some fixed metric on $M$, $f$ is Morse-Smale. One can also
weaken the question by asking $P$ and $Q$ to be fused to a "reasonable" 
critical point \cite{Cornea4}, a class that contains all critical points of locally analytic
functions. 

\begin{lem} Assume that $f^{-1}(-\infty,f(Q))$ has the homotopy type
of a $k$-dimensional $CW$-complex with $k< q-1$. If $H(P,Q)\not=0$, 
then $P$ and $Q$ cannot be fused to a reasonable critical point.
\end{lem}

\begin{proof}
This is a simple consequence of some results in \cite{Cornea4}.
First, as $P$ and $Q$ are consecutive and $f$ is Morse-Smale we may
assume that $f(P)>f(Q)$ and that $P$ and $Q$ are the only critical points in 
$f^{-1}[f(Q),f(P)]$. 
If $P$ and $Q$ can be fused to a reasonable critical point, one deduces \cite{Cornea4}
a cofibration sequence $Z\longrightarrow M'\longrightarrow M'''$
where $M'=f^{-1}(-\infty,f(Q)-\tau]$, $M'''=f^{-1}(-\infty, f(P)+\tau]$
with $\tau >0$  small. Let $M''=f^{-1}(-\infty, f(Q)+\tau]$.
 The composition $S^{p-1} \longrightarrow
M''\stackrel{\nabla}{\longrightarrow}S^{q}\vee M''$ is homotopic to
$S^{p-1}\stackrel{\nabla'}{\longrightarrow} S^{q}\vee
S^{p-1}\stackrel{id\vee\alpha (P)}{\longrightarrow}  S^{q}\vee M''$ where $\nabla'$
is the coaction associated to the cofibration sequence 
$S^{q-1}\longrightarrow Z\longrightarrow S^{p-1}$. As  $\pi_{p-1}(
S^{q}\vee S^{p-1})= \pi_{p-1}(S^{q}) \oplus \pi_{p-1}(S^{p-1})$ it follows that
$H(P,Q)$ vanishes. 
\end{proof}

\begin{rem} Controlling the behavior of Hopf
invariants has recently become a key tool in the homotopical study of the
Lusternik-Schnirelmann category. Results based on this technique are
the negative solution of the Ganea conjecture by Iwase \cite{Iwase}, the
examples of Roitberg showing that the L.S.-category is not generic in
the sense of the Mislin genus \cite{Roitberg} and the examples of Stanley  of spaces
of category $n$ but cone-length $n+1$ \cite{Stanley}. 
In all these examples the non-vanishing of certain Hopf invariants
associated to the attachment of a cell is used to deduce that the L.S.-category
is increased by the attachment of that cell. This is clearly also the homotopical
content of the lemma above.
\end{rem}

\bibliographystyle{amsplain}

\end{document}